\documentclass[11pt]{article}
\usepackage{}
\usepackage[total={6in, 9in}]{geometry}
 \usepackage{amsfonts}
\usepackage{amsthm}
\usepackage{amssymb}
\usepackage{amsmath,amsfonts}
\usepackage{mathrsfs}
\usepackage{cases}
\usepackage{latexsym,bm}
\usepackage{indentfirst}
\usepackage{color}
\usepackage{ifpdf}
\usepackage{graphicx}
\usepackage{psfrag}
\usepackage{dsfont}
\usepackage[
pdfauthor={},
pdftitle={},
pdfstartview=XYZ,
bookmarks=true,
colorlinks=true,
linkcolor=blue,
urlcolor=blue,
citecolor=blue,
bookmarks=true,
linktocpage=true,
hyperindex=true
]{hyperref}

\usepackage{graphicx}

\usepackage[natural]{xcolor}
\usepackage{tikz}
\usepackage{subfigure,tikz}
\usepackage{tikz-3dplot}
\usepackage{pgf,tikz,pgfplots}
\usepackage{mathrsfs}
\usepackage{enumerate} 

\newtheorem{THM}{\textbf{Theorem}}

\newtheorem{CLA}{\textbf{Claim}}

\newcommand{\pf}{\noindent\textbf{Proof}.\quad}

\newtheorem{FAC}{\textbf{Fact}}

\newcommand{\CC}{\mathcal{C}}

\usepackage{pdflscape}

\begin{document}
\title{A construction of a $\frac{3}{2}$-tough plane triangulation  with no 2-factor}
\author{Songling Shan \\ 
	\medskip  Auburn  University, Auburn, AL 36849\\
	\medskip 
	{\tt szs0398@auburn.edu}
}

\date{\today}
\maketitle

\emph{\textbf{Abstract}.}
In 1956, Tutte proved the celebrated theorem that every 4-connected planar graph is hamiltonian.
This result implies that every  more than $\frac{3}{2}$-tough planar graph
on at least three vertices is hamiltonian and so has a 2-factor.  Owens in 1999
constructed non-hamiltonian maximal planar graphs of 
toughness arbitrarily close to $\frac{3}{2}$ and asked  whether there exists a  maximal non-hamiltonian planar graph of toughness exactly $\frac{3}{2}$.  In fact, the graphs Owens constructed do not even contain a 2-factor. Thus the toughness of exactly  $\frac{3}{2}$  is the only case left in asking the existence of 2-factors in  tough planar graphs. This question was also asked by  Bauer, Broersma,
and  Schmeichel in a survey.  In this paper, we 
close this gap by constructing a  maximal  $\frac{3}{2}$-tough plane graph  with no  2-factor, 
 answering  the  question asked by Owens as well as by  Bauer, Broersma,
 and  Schmeichel.


\vspace{2mm}

\section{Introduction}

We consider only simple  graphs. 
%The notation defined is for both simple graphs 
%and multigraphs although we will not emphasis 
%on that.  
Let $G$ be a graph.
Denote by $V(G)$ and  $E(G)$ the vertex set and edge set of $G$,
respectively. We denote by $n(G)$ and $e(G)$ the sizes of 
$V(G)$ and $E(G)$, respectively, and by $f(G)$
the number of faces of $G$ if $G$ is embedded on a surface.  
Let $v\in V(G)$ and $S\subseteq V(G)$. 
Then  $N_G(v)$   denotes the set of neighbors
of $v$ in $G$ and  $N_G(S):=(\bigcup_{x\in S}N_G(x))\setminus S$. 
The subgraph of $G$ induced on $S$ and $V(G)\setminus S$ are denoted by
 $G[S]$ and $G-S$, respectively. For notational simplicity we write $G-x$ for $G-\{x\}$.
Let $V_1,
V_2\subseteq V(G)$ be two disjoint vertex sets. Then $E_G(V_1,V_2)$ is the set
of edges in $G$  with one end in $V_1$ and the other end in $V_2$ and  $e_G(V_1,V_2):=|E_G(V_1,V_2)|$.  We write $E_G(v,V_2)$ and $e_G(v,V_2)$
if $V_1=\{v\}$ is a singleton.  
%We also use $G[V_1,V_2]$
%to denote the bipartite subgraph of $G$ with vertex set $V_1\cup V_2$
%and edge set $E_G(V_1, V_2)$. 
When $H\subseteq G$ and $S\subseteq V(G)\setminus V(H)$, we write $N_G(H)$, $E_G(H,S)$,
and $e_G(H,S)$ respectively for $N_G(V(H))$, $E_G(V(H),S)$,
and $e_G(V(H),S)$.  For two integers $p$ and $q$, we let $[p,q]=\{i\in \mathbb{Z}:  p\le i\le q\}$.

The number of components of $G$ is denoted by $c(G)$. Let $t\ge 0$ be a
real number. The graph $G$ is said to be \emph{$t$-tough} if $|S|\ge t\cdot
c(G-S)$ for each $S\subseteq V(G)$ with $c(G-S)\ge 2$. The \emph{toughness $\tau(G)$} is the largest real number $t$ for which $G$ is
$t$-tough, or is  $\infty$ if $G$ is complete. This concept was introduced by Chv\'atal~\cite{chvatal-tough-c} in 1973.
 It is  easy to see that  if $G$ has a hamiltonian cycle
then $G$ is 1-tough. Conversely,
 Chv\'atal~\cite{chvatal-tough-c}
 conjectured that
there exists a constant $t_0$ such that every
$t_0$-tough graph is hamiltonian.
 Bauer, Broersma and Veldman~\cite{Tough-counterE} have constructed
$t$-tough graphs that are not hamiltonian for all $t < \frac{9}{4}$, so
$t_0$ must be at least $\frac{9}{4}$ if Chv\'atal's conjecture is true.
The conjecture has
been verified when restricted to a number of graph
classes~\cite{Bauer2006},
including planar graphs, claw-free graphs, co-comparability graphs, and
chordal graphs.   

The study of cycle structures 
in planar graphs under a given toughness condition is particularly intensive and 
interesting, see, for examples~\cite{ MR1115734, MR1248489,MR1364524,  MR1724905,  MR1395475}.   Observe that  any more than $\frac{3}{2}$-tough planar graph on at least 5 vertices 
is 4-connected. Thus the toughness conjecture of  Chv\'atal holds 
for planar graphs with  toughness   greater than $\frac{3}{2}$ by the classic result of 
Tutte~\cite{MR81471} that every 4-connected planar graph is hamiltonian. 
Furthermore,  it is shown by Owens~\cite{MR1724905} that $t_0$ cannot be smaller 
than $\frac{3}{2}$. It is still unknown whether $t_0=\frac{3}{2}$
is the sharp toughness bound for a planar graph to be hamiltonian. 
In fact, this question is even open for the existence of 2-factors in planar graphs. 

A \emph{2-factor} in a graph $G$ is a spanning 2-regular subgraph. 
Thus a  hamiltonian cycle of $G$ is a 2-factor with only one component.  
By the result of Tutte~\cite{MR81471}, we know that  every more than $\frac{3}{2}$-tough planar graph on at least 3 vertices has a 2-factor. On the other hand, constructed by Owens~\cite{MR1724905}, there are maximal planar graphs with toughness arbitrarily close to $\frac{3}{2}$ but with no  2-factor.  Owens  asked  in the same paper whether there exists a non-hamiltonian maximal planar graph with toughness exactly $\frac{3}{2}$. 
Bauer,  Broersma, and Schmeichel in the survey~\cite{Bauer2006} 
commented that 
``one of the challenging open problems in this area is to determine whether
every $\frac{3}{2}$-tough maximal planar graph has a 2-factor. If so, are they all hamiltonian? We also
do not know if a $\frac{3}{2}$-tough planar graph has a 2-factor.'' 
In this paper,  we answer positively the question asked by Owens and negatively the questions raised by Bauer,  Broersma, and Schmeichel. 
%The result provides in a sense also some evidence for the existence of a hamiltonian cycle in a 
%$\frac{3}{2}$-tough maximal planar graph.   

\begin{THM}\label{main}
There exists a $\frac{3}{2}$-tough plane triangulation with no  2-factor. 
\end{THM}

The remainder of this paper is organized as follows: in Section 2, we introduce some notation and preliminary 
results, and in Section 3, we prove Theorem~\ref{main}. 

\section{Preliminary results}

%In this section,  we introduce Tutte's 2-factor Theorem and a result on matchings. 

%\subsection{Tutte's 2-Factor Theorem}
Let $S$ and $T$ be disjoint sets of vertices of a graph $G$ 
and  $D$ be a component of $G-(S\cup T)$.
Then $D$ is said to be an \emph{odd component}
(resp.~\emph{even component})
if $e_G(D, T)\equiv 1\pmod{2}$
(resp.~$e_G(D, T)\equiv 0\pmod{2}$).
For each integer $k\ge 0$, we denote by $\CC_{2k+1}$ the set of all odd components $D$ of 
$G-(S\cup T)$ such that $e_G(D,T)=2k+1$. Let $\CC=\bigcup_{k\ge 0} \CC_{2k+1}$ and 
 let $c(S, T)=|\CC|$.

Define 
$\delta(S, T)=2|S|+\sum_{y\in T} d_{G-S}(y)-2|T|-c(S, T)$.
It is easy to see $\delta(S, T)\equiv 0\pmod{2}$
for every $S$,~$T\subseteq V(G)$
with $S\cap T=\emptyset$.
We use the following criterion for the existence of a $2$-factor,
which is a special case of Tutte's $f$-Factor Theorem.
\begin{THM}[Tutte~\cite{Tutte}]\label{tutte's theorem}
	A graph $G$ has a $2$-factor if and only if
	$\delta(S, T)\ge 0$
	for every $S$,~$T\subseteq V(G)$
	with $S\cap T=\emptyset$.
\end{THM}

An ordered pair $(S,T)$ consisting of disjoint sets of vertices $S$ and $T$ in a graph $G$ 
is called a barrier if $\delta_G(S, T)\le -2$.
By Theorem~\ref{tutte's theorem},
if $G$ does not have a $2$-factor,
then $G$ has a barrier.
%A barrier $(S, T)$ is called a \emph{biased barrier}
%if $|T|$ is  minimum and  subject to that $|S|$ is maximum among all the barriers of $G$.
%We will use the following properties of a biased barrier.  (The 
%property was proved in~\cite{MR3956455}  
%for a barrier $(S,T)$ defined as $|S|$ is  maximum and  subject to that $|T|$ is minimum among all the barriers of $G$, but the proof works 
%for both definitions.) 
%
%
%\begin{LEM}[Lemma 3.2, \cite{MR3956455}]\label{lem:minimal_barrier}
%	Let $G$ be a graph with no  $2$-factor,
%	and let $(S, T)$ be a biased barrier of $G$.
%	Then
%	\begin{enumerate}
%		\item
%		$T$ is independent in $G$,
%		\item if $D$ is an even component
%		with respect to $(S, T)$,
%		then $e_G(T, D)=0$,
%		\item
%		if $D$ is an odd component
%		with respect to $(S, T)$,
%		then $e_G(y, D) \le 1$
%		for every $y\in T$,
%		\item
%		if $D$ is an odd component
%		with respect to $(S, T)$,
%		then $e_G(x, T) \le 1$
%		for every $x\in V(D)$.  
%	\end{enumerate}
%\end{LEM}

We need also the following result regarding the toughness of the square of a graph. 
%\begin{THM}[M.~B.~Dillencourt, {\cite[Lemma 2]{MR1248489}}]\label{lem:component-in-triangulation}
%Let $G$ be a plane triangulation and $S\subseteq V(G)$.  Then for any face $F$ of $G[S]$, the subgraph 
%of $G$ induced by those vertices of $G$ lying inside $F$ is a connected graph. 
%\end{THM}

\begin{THM}[V.~Chv\'atal, {\cite[Theorem 1.7]{chvatal-tough-c}}]\label{lem:toughness-square}
For any graph $G$, we have $\tau(G^2) \ge \kappa(G)$, where $G^2$ is obtained from $G$ by adding edges 
joining pairs of vertices of distance 2 in $G$, and $\kappa(G)$ is the connectivity of $G$.  
\end{THM}

\section{Proof of Theorem~\ref{main}}

Suppose there exists a $\frac{3}{2}$-tough plane triangulation on $n$ vertices with no 2-factor. 
Let $(S,T)$ be a barrier of $G$, and $\CC$ be the set of all odd components of $G-(S\cup T)$. 
Let $G_2$ be obtained from $G$ by deleting all the vertices in $S$;
$G_1$ be obtained from $G_2$ by smoothing  all the  vertices of $T$ that have degree 2 in $G_2$  (for a degree 2 vertex $u$, smooth it amounts to deleting $u$  and joining an edge between its two neighbors) and deleting all the  vertices of $T$ that have degree 1 in $G_2$;  and $G_0$  be  obtained from $G_1$ by contracting each graph  $D$ in $\CC$ into a single vertex (identify all the vertices of $D$ into a single vertex and join all the edges from the vertex to  vertices of $N_{G_1}(D)$). 
We call $G_0$  the \emph{component graph} of $G$.  For any subgraph $G'$ of $G$,  
we also call 
the subgraph of $G_0$ obtained through the above process from $G'$ (replacing  ``contracting each graph in $\CC$'' by  ``contracting the restriction of  each graph of  $\CC$ in $G'$'') the \emph{component graph of $G'$}. 
  Our construction of $G$ starts with  
its component graph $G_0$ and then  reverse the process of getting $G_0$ from $G$ to obtain the graph $G$.  In the construction, we will add vertices
 and edges to the existing graph in a way such that  the resulting graph is still a plane graph.  We will stick to this rule without mentioning it. 

{\noindent {\bf Step 1}: The component graph $G_0$.}

\medskip

\begin{landscape}
	
	\begin{figure}[!htb]
		\vspace{-3cm}
		\begin{center}
			\begin{tikzpicture}[scale=1]	
				\usetikzlibrary{calc}
				\begin{scope}[shift={(0,0)}, rotate=0]
					
					%define node coordinates 
					
					\def \n {86}
					\def \radius {7cm}
					\def \margin {4}
					\foreach \s in {1,...,\n}
					{
						{\tikzstyle{every node}=[draw ,circle,fill=white, minimum size=0.35cm,
							inner sep=1pt]
							\node[draw] at ({360/\n *\s - 1)}:\radius)(\s) {};
							\draw[-, >=latex] ({360/\n * (\s -1)+1.6\margin}:\radius) 
							arc ({360/\n * (\s - 1)+0.55*\margin}:{360/\n * (\s)-0.8*\margin}:\radius);
						}
					}

					%drawing vertices
					
					\foreach \s in {1,3,5,7,9,11,13,15,17,19,21,23,25,27,29,31,33,35,37,39,41,43,45,47,49,51,53,55,57,59,61,63,65,67,69,71,73,75,77,79,81,83,85}
					{\node [circle, fill=black!20!white, label={} ]at (\s) {};
					}
					
					\foreach \s in {5,9,13,17, 25,29,33,37,45,49,53,57,65, 69,73,77,81}
					{\node [circle, fill=black, label={} ]at (\s) {};
					}
					
					\foreach \s in {51,55,59,75,79, 83, 11,15,19,31,35,39}
					{\node [circle, fill=gray, label={} ]at (\s) {};
					}
					
					%\node [circle, fill=black, label={[label distance=-0.1cm]0:$v_1$}] at (1) {};
					%\node [circle,  label={[label distance=-0.1cm]0:$v_{86}$}] at (86) {};
					%\node [circle,  label={[label distance=-0.1cm]right:$v_2$}] at (2) {};
					%\node [circle,  label={[label distance=-0.1cm]0:$v_{3}$}] at (3) {};
					%\node [circle,  label={[label distance=-0.1cm]0:$v_{4}$}] at (4) {};
					%\node [circle,  label={[label distance=-0.1cm]0:$v_{5}$}] at (5) {};
					%\node [circle,  label={[label distance=-0.1cm]0:$v_{6}$}] at (6) {};
					%\node [circle,  label={[label distance=-0.1cm]0:$v_{7}$}] at (7) {};
					%\node [circle,  label={[label distance=-0.1cm]0:$v_{8}$}] at (8) {};
					%\node [circle,  label={[label distance=-0.1cm]0:$v_{9}$}] at (9) {};
					%\node [circle,  label={[label distance=-0.1cm]0:$v_{10}$}] at (10) {};
					
					\foreach \s in {1}{
						\node [circle,  label={[label distance=0.07cm]0:$v_{{\s}}$}] at (\s) {};
					}
					
					\foreach \s in {3, 5, 7,  11}{
						\node [circle,  label={[label distance=-0.05cm]0:$v_{{\s}}$}] at (\s) {};
					}
					
					\foreach \s in {9}{
						\node [circle,  label={[label distance=-0.22cm]20:$v_{{\s}}$}] at (\s) {};
					}
					
					\foreach \s in {13}{
						\node [circle,  label={[label distance=-0.15cm]45:$v_{{\s}}$}] at (\s) {};
					}

					\foreach \s in {15,17, 21}{
						\node [circle,  label={[label distance=-0.1cm]90:$v_{{\s}}$}] at (\s) {};
					}
					
					\foreach \s in {19}{
						\node [circle,  label={[label distance=0.02cm]90:$v_{{\s}}$}] at (\s) {};
					}
					
					\foreach \s in {23}{
						\node [circle,  label={[label distance=0cm]90:$v_{{\s}}$}] at (\s) {};
					}
					
					\foreach \s in {25,27,29}{
						\node [circle,  label={[label distance=-0.1cm]90:$v_{{\s}}$}] at (\s) {};
					}
					
					\foreach \s in {31,35}{
						\node [circle,  label={[label distance=-0.1cm]150:$v_{{\s}}$}] at (\s) {};
					}
					
					\foreach \s in {37}{
						\node [circle,  label={[label distance=-0.15cm]150:$v_{{\s}}$}] at (\s) {};
					}
					
					\foreach \s in {33}{
						\node [circle,  label={[label distance=-0.15cm]140:$v_{{\s}}$}] at (\s) {};
					}
					
					\foreach \s in {39}{
						\node [circle,  label={[label distance=0cm]180:$v_{{\s}}$}] at (\s) {};
					}
					
					\foreach \s in {41,43,45,49}{
						\node [circle,  label={[label distance=-0.1cm]180:$v_{{\s}}$}] at (\s) {};
					}

					\foreach \s in {47}{
						\node [circle,  label={[label distance=0.02cm]178:$v_{{\s}}$}] at (\s) {};
					}
					
					\foreach \s in {51,53,55,57,57}{
						\node [circle,  label={[label distance=-0.1cm]200:$v_{{\s}}$}] at (\s) {};
					}
					
					\foreach \s in {59, 61,63,65,67}{
						\node [circle,  label={[label distance=-0.1cm]270:$v_{{\s}}$}] at (\s) {};
					}
					
					\foreach \s in {69,71,73}{
						\node [circle,  label={[label distance=-0.1cm]90:$v_{{\s}}$}] at (\s) {};
					}
					
					\foreach \s in {75,77,79,81,83}{
						\node [circle,  label={[label distance=-0.1cm]100:$v_{{\s}}$}] at (\s) {};
					}
					
					\foreach \s in {85}{
						\node [circle,  label={[label distance=-0.1cm]0:$v_{{\s}}$}] at (\s) {};
					}
					
					%\coordinate (d) at (0,0);
					%\node[circle, fill=white, label={above:$w$}] at (0,0) {};
					
					\node[draw, circle,fill=white] at (0,0)(d) {$w$}; 
					
					%	\foreach \s in {1,...,85}
					%	{
						%		\pgfmathparse{\s+1};
						%		\xdef\t{\pgfmathresult};
						%		\path[draw,black]
						%	(\s) edge node[name=la,pos=0.7, above] {\color{blue} } (\t); 
						%}
					
					\path[draw,black]
					(d) edge node[name=la,pos=0.7, above] {\color{blue} } (1)
					(d) edge node[name=la,pos=0.7, above] {\color{blue} } (4)
					(d) edge node[name=la,pos=0.7, above] {\color{blue} } (6)
					(d) edge node[name=la,pos=0.7, above] {\color{blue} } (8)
					(d) edge node[name=la,pos=0.7, above] {\color{blue} } (10)
					(d) edge node[name=la,pos=0.7, above] {\color{blue} } (12)
					(d) edge node[name=la,pos=0.7, above] {\color{blue} } (14)
					(d) edge node[name=la,pos=0.7, above] {\color{blue} } (16)
					(d) edge node[name=la,pos=0.7, above] {\color{blue} } (18)
					;

					{\tikzstyle{every node}=[draw ,circle,fill=black, minimum size=0.3cm,
						inner sep=0pt]
						\draw [] (d) -- (4) node [midway] (d1){};
						\draw [] (d) -- (6) node [midway] (d2){};
						\draw [] (d) -- (8) node [midway] (d3){};
						\draw [] (d) -- (10) node [midway] (d4){};
						\draw [] (d) -- (12) node [midway](d5) {};
						\draw [] (d) -- (14) node [midway] (d6){};
						\draw [] (d) -- (16) node [midway] (d7){};
						\draw [] (d) -- (18) node [midway] (d8){};
					}

					\foreach \s in {1}{
						\node [circle,  label={[label distance=-0.1cm]-90:$u_{{\s}}$}] at (d\s) {};
					}
					
					\foreach \s in {2,3,4,5}{
						\node [circle,  label={[label distance=-0.09cm]-1:$u_{{\s}}$}] at (d\s) {};
					}
					
					\foreach \s in {6,7}{
						\node [circle,  label={[label distance=-0.1cm]3:$u_{{\s}}$}] at (d\s) {};
					}

					\foreach \s in {8}{
						\node [circle,  label={[label distance=-0.1cm]49:$u_{{\s}}$}] at (d\s) {};
					}
					
					\path[draw,black]
					(d) edge node[name=la,pos=0.7, above] {\color{blue} } (21)
					(d) edge node[name=la,pos=0.7, above] {\color{blue} } (24)
					(d) edge node[name=la,pos=0.7, above] {\color{blue} } (26)
					(d) edge node[name=la,pos=0.7, above] {\color{blue} } (28)
					(d) edge node[name=la,pos=0.7, above] {\color{blue} } (30)
					(d) edge node[name=la,pos=0.7, above] {\color{blue} } (32)
					(d) edge node[name=la,pos=0.7, above] {\color{blue} } (34)
					(d) edge node[name=la,pos=0.7, above] {\color{blue} } (36)
					(d) edge node[name=la,pos=0.7, above] {\color{blue} } (38)
					
					;

					{\tikzstyle{every node}=[draw ,circle,fill=black, minimum size=0.3cm,
						inner sep=0pt]
						\draw [] (d) -- (24) node [midway](d9) {};
						\draw [] (d) -- (26) node [midway] (d10){};
						\draw [] (d) -- (28) node [midway] (d11){};
						\draw [] (d) -- (30) node [midway] (d12){};
						\draw [] (d) -- (32) node [midway] (d13){};
						\draw [] (d) -- (34) node [midway] (d14){};
						\draw [] (d) -- (36) node [midway] (d15){};
						\draw [] (d) -- (38) node [midway] (d16){};
						
					}

					\foreach \s in {9,10}{
						\node [circle,  label={[label distance=-0.1cm]95:$u_{{\s}}$}] at (d\s) {};
					}
					
					\foreach \s in {11,12}{
						\node [circle,  label={[label distance=-0.1cm]170:$u_{{\s}}$}] at (d\s) {};
					}

					\foreach \s in {13,14,15,16}{
						\node [circle,  label={[label distance=-0.1cm]180:$u_{{\s}}$}] at (d\s) {};
					}
					
					\path[draw,black]
					(d) edge node[name=la,pos=0.7, above] {\color{blue} } (41)
					(d) edge node[name=la,pos=0.7, above] {\color{blue} } (44)
					(d) edge node[name=la,pos=0.7, above] {\color{blue} } (46)
					(d) edge node[name=la,pos=0.7, above] {\color{blue} } (48)
					(d) edge node[name=la,pos=0.7, above] {\color{blue} } (50)
					(d) edge node[name=la,pos=0.7, above] {\color{blue} } (52)
					(d) edge node[name=la,pos=0.7, above] {\color{blue} } (54)
					(d) edge node[name=la,pos=0.7, above] {\color{blue} } (56)
					(d) edge node[name=la,pos=0.7, above] {\color{blue} } (58)
					
					;

					{\tikzstyle{every node}=[draw ,circle,fill=black, minimum size=0.3cm,
						inner sep=0pt]
						\draw [] (d) -- (44) node [midway] (d17){};
						\draw [] (d) -- (46) node [midway] (d18){};
						\draw [] (d) -- (48) node [midway] (d19){};
						\draw [] (d) -- (50) node [midway] (d20){};
						\draw [] (d) -- (52) node [midway] (d21){};
						\draw [] (d) -- (54) node [midway] (d22){};
						\draw [] (d) -- (56) node [midway](d23) {};
						\draw [] (d) -- (58) node [midway] (d24){};
						
					}

					\foreach \s in {17,18}{
						\node [circle,  label={[label distance=-0.15cm]175:$u_{{\s}}$}] at (d\s) {};
					}
					
					\foreach \s in {19}{
						\node [circle,  label={[label distance=-0.05cm]179:$u_{{\s}}$}] at (d\s) {};
					}
					
					\foreach \s in {20,21,22}{
						\node [circle,  label={[label distance=-0.15cm]180:$u_{{\s}}$}] at (d\s) {};
					}

					\foreach \s in {23}{
						\node [circle,  label={[label distance=0cm]270:$u_{{\s}}$}] at (d\s) {};
					}
					
					\foreach \s in {24}{
						\node [circle,  label={[label distance=-0.1cm]0:$u_{{\s}}$}] at (d\s) {};
					}

					\path[draw,black]
					(d) edge node[name=la,pos=0.7, above] {\color{blue} } (61)
					(d) edge node[name=la,pos=0.7, above] {\color{blue} } (64)
					(d) edge node[name=la,pos=0.7, above] {\color{blue} } (66)
					(d) edge node[name=la,pos=0.7, above] {\color{blue} } (68)
					(d) edge node[name=la,pos=0.7, above] {\color{blue} } (70)
					(d) edge node[name=la,pos=0.7, above] {\color{blue} } (72)
					(d) edge node[name=la,pos=0.7, above] {\color{blue} } (74)
					(d) edge node[name=la,pos=0.7, above] {\color{blue} } (76)
					(d) edge node[name=la,pos=0.7, above] {\color{blue} } (78)
					(d) edge node[name=la,pos=0.7, above] {\color{blue} } (80)
					(d) edge node[name=la,pos=0.7, above] {\color{blue} } (82)
					(d) edge node[name=la,pos=0.7, above] {\color{blue} } (85)
					;

					{\tikzstyle{every node}=[draw ,circle,fill=black, minimum size=0.3cm,
						inner sep=0pt]
						\draw [] (d) -- (64) node [midway] (d25){};
						\draw [] (d) -- (66) node [midway] (d26){};
						\draw [] (d) -- (68) node [midway] (d27){};
						\draw [] (d) -- (70) node [midway] (d28){};
						\draw [] (d) -- (72) node [midway] (d29){};
						\draw [] (d) -- (74) node [midway] (d30){};
						\draw [] (d) -- (76) node [midway](d31) {};
						\draw [] (d) -- (78) node [midway] (d32){};
						\draw [] (d) -- (80) node [midway] (d33){};
						\draw [] (d) -- (82) node [midway] (d34){};
						
					}

					\foreach \s in {25}{
						\node [circle,  label={[label distance=-0.1cm]180:$u_{{\s}}$}] at (d\s) {};
					}

					\foreach \s in {26}{
						\node [circle,  label={[label distance=-0.25cm]-30: $u_{{\s}}$}] at (d\s) {};
					}
					
					\foreach \s in {27}{
						\node [circle,  label={[label distance=-0.2cm]-30: $u_{{\s}}$}] at (d\s) {};
					}
					
					\foreach \s in {28,29,30}{
						\node [circle,  label={[label distance=-0.1cm]-10: $u_{{\s}}$}] at (d\s) {};
					}

					\foreach \s in {31,32,33,34}{
						\node [circle,  label={[label distance=-0.1cm]0: $u_{{\s}}$}] at (d\s) {};
					}
					
					%\path[draw,black]
					%(42) to [out=-60,in=-90] (55)
					%(40.south) to [out=270,in=200] (3,-8.5)
					%(3,-8.5) to [out=20,in=290] (59)
					%; 
					%
					%\path[draw,black]
					%(35) to [out=150,in=180] (22)
					%(39.west) to [out=180,in=240] (-7.5,4.5)
					%(-7.5,4.5) to [out=60,in=180] (20)
					%; 
					%
					%%\node at (21,1.4) {$c$}; 
					%%\node at (21,-3.4) {$d$}; 
					%%
					%%\node at (7.2,1) {$A$};
					%%\node at (7.2,-3) {$B$};
					%
					%\path[draw,black]
					%(15) to [out=0,in=150] (5.6, 6)
					%(5.6,6) to [out=-30,in=100]  (7.8,1.6)
					%(7.8, 1.6) to [out=280,in=0] (2)
					%(19) to [out=20,in=120] (7.3,5.5)
					%(7.3,5.5) to[out=-60,in=60] (8, 1)
					%(8,1) to [out=240,in=30] (86)
					%; 
					
					\draw[bend right=59,-]  (7) to node [auto] {} (20);
					\draw[bend right=100,-]  (3) to node [auto] {} (22);
					\draw[bend right=60,-]  (27) to node [auto] {} (40);
					\draw[bend right=110,-]  (23) to node [auto] {} (42);
					
					\draw[bend right=125,-]  (47) to node [auto] {} (60);
					\draw[bend right=110,-]  (43) to node [auto] {} (62);
					
					\draw[bend right=55,-]  (67) to node [auto] {} (86);
					\draw[bend right=90,-]  (63) to node [auto] {} (2);
					\draw[bend right=40,-]  (71) to node [auto] {} (84);
					
					%\draw [red] plot [smooth] coordinates {(15) (5.6,6) (7.8,1.6) (2)};
					
					%\draw [-] (-360/86:7.1)  arc (-360/86:360/15:9.1) node [right,pos=0.5] {};
					%\draw (86) arc 
					%[
					%start angle=30,
					%end angle=80,
					%x radius=0.2cm,
					%y radius =10cm
					%] ;
					
					\thispagestyle{empty}
				\end{scope}	
			\end{tikzpicture}
		\end{center}
		\vspace{-2cm}
		\caption{The component graph $G_0$ of $G$.}
		\label{fig1}
	\end{figure}
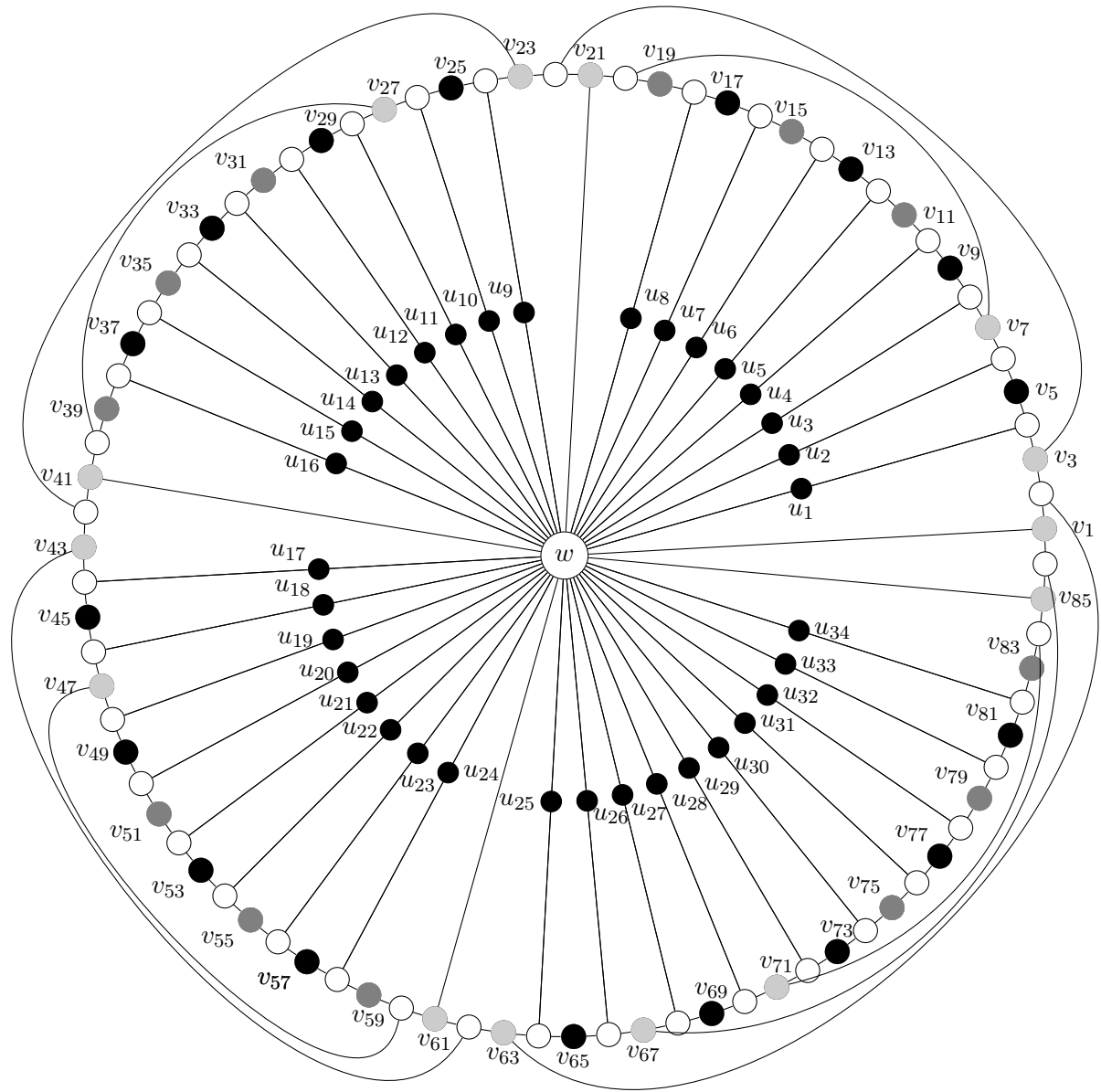
\end{landscape}

Let  $G_0$ be  the plane graph drawn in Figure~\ref{fig1}.  
The vertices of $G_0$ are depicted in 4 different colors: white vertices (those vertices form a color class of $G_0$), light gray vertices, gray vertices, and black vertices, where all vertices of 
gray or  black color have degree 2 in $G_0$.  The  filled vertices are labelled while the white vertices (except $w$) are not labelled due to the space limitation.  However, we will  denote  the  white vertices 
between two filled vertices $v_i$ and $v_{i+2}$ by $v_{i+1}$ for each $i\in [1, 85]$, where $v_{87} :=v_1$.  We denote by $C_0$ the cycle $v_1v_2\ldots v_{86}v_1$. 
It is easy to see that the graph is bipartite 
with the bipartition as the set of filled vertices and the set of unfilled vertices in the drawing.  
As  there are only 43 unfilled vertices in $G_0-w$, we have 
 the following fact.  

\begin{FAC}\label{fact:matching-size-in-G0}
The size of a maximum matching in $G_0-w$ is 43. 
\end{FAC}

Each face of $G_0$ that has  exactly three degree 2 vertices of the same color on its boundary is called an \emph{$S$-triangle face} associated with those degree 2 vertices. For  example, the face with boundary  $wu_1v_4v_5v_6u_2w$ is associated with $u_1, u_2, v_5$,
and the face with boundary $v_7v_8\ldots v_{20}v_7$ is associated with  $v_{11}, v_{15}, v_{19}$.

The vertex $w$ in $G_0$ is corresponding to the  only component in $\CC_{39}$, 
and all other vertices of $G_0$ are corresponding to  graphs in $\CC_3$, which are all 
triangles. Thus,  we can get to our second step of construction by replacing each vertex of 
$G_0$ with a graph.

\medskip 

{\noindent {\bf Step 2}: The  construction of  $G_1$.}

\medskip 

The vertex $w$ will be replaced with a 4-connected plane graph. 
Let a plane graph $D$ be constructed as follows:
\begin{enumerate}[(i)]
	\item Take  four  vertex-disjoint 39-cycles  $A_1, A_2, A_3, A_4$ with vertex set $V(A_i)=\{a_{i,1}, \ldots, a_{i,39}\}$ 
	for $i\in [1,4]$ 
	such that they are embedded in   the plane with the faces  bounded by $A_1$, $A_1$ and $A_2$, $A_2$ and $A_3$,  $A_3$ and $A_4$, and $A_4$, respectively; 
	\item Add edges $a_{i,j}a_{(i+1),j}$ for each $i\in [1,3]$ and  $j\in [1,39]$; 
	\item Add edges $a_{i,j}a_{(i+1),(j+1)}$ for each $i\in [1,3]$ and  $j\in [1,39]$, where $a_{(i+1),40}:=a_{(i+1),1}$; 
	\item Add edges  $a_{4,1} a_{4,j}$ for each $j\in [3,38]$. 
\end{enumerate}
Note that   $D$ is a near triangulation with $A_1$  being   the boundary of its non-triangle face.  
We prove that $D$ is 2-tough. This fact will be used later on when we show that the finally constructed graph $G$ is at least $\frac{3}{2}$-tough. 

\begin{CLA}\label{fact:D-2-tough}
	The graph	$D$ is 2-tough.  
\end{CLA}

\pf   
%We first show that $D$ is 4-connected. 
%The graph $D$ is clearly 2-connected with $\delta(D)=4$. 
%We verify that $D$ has no 3-cut. Suppose otherwise that $D$ has a cut $W$ of size 3. 
%As $\delta(G)=4$, each component of $D-W$ contains an edge.  As each vertex of $A\cup C$ has in $D$ two  
%neighbors  from $V(B)$, it follows that $|W\cap V(B)| \ge 2$.
%Thus  $|W\cap (V(A) \cup V(C))| \le 1$. This implies that both $A-W$ and $C-W$ are connected, 
%and each vertex of $V(B)\setminus W$ has a neighbor from each of $V(A)\setminus W$ and $V(C)\setminus W$ in $D-W$. 
%Therefore $D-W$ is connected, a contradiction. Thus $D$ is 4-connected. 
Let $D_1=D[V(A_1)\cup V(A_2)]$ and $D_2=D[V(A_3)\cup V(A_4)]$. 
We first show that each $D_i$ for $i\in [1,2]$ is 2-tough. As $D_2$ contains a spanning subgraph isomorphic to $D_1$ by the construction of $D$ and so $\tau(D_2) \ge \tau(D_1)$, 
we only need to show that $\tau(D_1) \ge 2$. By Theorem~\ref{lem:toughness-square}, we show that $D_1$ is the square of a 78-cycle. 
This is obvious as if we  let $$Q=a_{1,1} a_{2,2}a_{1,2}a_{2,3} \ldots a_{1,i} a_{2,(i+1)} a_{1,(i+1)} \ldots a_{1,38} a_{2,39} a_{1,39} a_{2,1} a_{1,1}$$  
be a 78-cycle, then  we can easily check that any two vertices of distance 2 in $Q$ are adjacent in $D_1$, and the edges of $G_1$ consists of the edges of $Q$,
and edges joining pairs of vertices that are of distance 2 in $Q$.  Thus $D_1 =Q^2$. 

Now we have $\tau(D_2) \ge \tau(D_1) \ge 2$,  and we show that $\tau(D) \ge 2$. Let $W$ be an arbitrary  cutset of $D$.  
Let $W_i=W\cap V(D_i)$ for $i\in [1,2]$. If $W_i$ is a cutset of $D_i$ for each $i\in [1,2]$, then we have $ |W_i| \ge 2 c(D_i-W_i) $. 
As $c(D-W) \le c(D_1-W_1)+c(D_2-W_2)$, it follows that $|W| \ge 2c(D-W)$.   

Consider next that $W_i$ is not  a cutset of $D_i$ for each $i\in [1,2]$. Then   $c(D_i-W_i) =1$,  and so $c(D-W) =2$ as $W$ is 
a cutset of $D$. If $|W| \ge 4$, then we have $|W| \ge 2c(D-W)$ already. Thus we assume  $|W| \le 3$, and assume by symmetry, that $|W_1| \le 1$. 
As $D_2-W_2$ is a graph that contains at least $39-3=36$ vertices from $V(A_3)$, and $D_1-W_1$ contains at least $39-1$ vertices from $V(A_2)$, 
it follows that $E_{D-W}(A_2-W_1, A_3-W_2) \ne \emptyset$ by the construction of $D$. Thus  $E_{D-W}(D_1-W_1, D_2-W_2) \ne \emptyset$ and so $c(D-W) =1$, a contradiction to the assumption that $W$ is a cutset of $D$. 

Consider lastly, by symmetry, that $W_1$ is a cutset of $D_1$ but $W_2$ is not a cutset of $D_2$.  Thus $c(D_2-W_2) =1$. 
We may further assume  $|W_2| \le 1$.  For otherwise,  we have $ |W| =|W_1|+|W_2|  \ge 2(c(D_1-W_1)+1 ) \ge  2c(D-W)$ already.  
As each vertex  from $V(A_2)$ has in $D$ two neighbors from $V(A_3)$, if $D_1-W_1$ has a component containing a vertex from $V(A_2)$, 
then we have $c(D-W)=c(D_1-W_1)$ and so  we get $|W| \ge 2c(D-W)$. Thus we assume that every component of $D_1-W_1$ is disjoint with $A_2$. 
As a consequence, $V(A_2) \subseteq W$. As $A_1$ is a cycle and so is 1-tough, we know that $|W_1\cap V(A_1)| \ge c(A_1-(W_1\cap V(A_1)))$. 
Thus $c(D-W) =c(A_1-(W_1\cap V(A_1)))+1$. As $c(A_1-(W_1\cap V(A_1))) \ge 2$ by $W_1$ being a cutset of $D_1$, it follows that  $|W_1\cap V(A_1)|  \le 37$. 
Hence  $|W| \ge 39+|W_1\cap V(A_1)| \ge 2(c(A_1-(W_1\cap V(A_1)))+1) =2c(D-W)$, as desired. 
\qed

%$D[V(A\cup B)]$ is the square of $C_{54}$, and $D[V(B\cup C)]$ is the square of $C_{54}$. 

{\noindent {\bf Step 2.1}: In $G_0$, we replace $w$ with $D$. That is, delete $w$,  place $D$ in the position of $w$ such that $G_0-w$ is embedded inside the non-triangle face of $D$, then add a perfect matching between $V(A)$
	and $N_{G_0}(w)$ such that the resulting graph is still a plane graph. }

\begin{figure}[!htb]
		\begin{center}
			\begin{tikzpicture}[scale=1]	
				\usetikzlibrary{calc}
				\begin{scope}[shift={(0,0)}, rotate=0]
				
						{\tikzstyle{every node}=[draw ,circle,fill=white, minimum size=0.35cm, inner sep=1pt]						
					\node[draw] at (0,0) (w){$w$};
					\node[draw] at (-1,-2) (u1){$u_1$};
					\node[draw] at (1,-2) (u2){$u_2$};
					\node[draw] at (-2,-4) (v4){$v_4$};
					\node[draw] at (2,-4) (v6){$v_6$};
						\node[draw] at (0,-4) (v5){$v_5$};
					}
					
		\path[draw,black]
	(w) edge node[name=la,pos=0.7, above] {\color{blue} } (u1)				
	(w) edge node[name=la,pos=0.7, above] {\color{blue} } (u2)		
	(u1) edge node[name=la,pos=0.7, above] {\color{blue} } (v4)		
	(u2) edge node[name=la,pos=0.7, above] {\color{blue} } (v6)		
	(v5) edge node[name=la,pos=0.7, above] {\color{blue} } (v6)		
	(v5) edge node[name=la,pos=0.7, above] {\color{blue} } (v4);

	\node[] at (0,-4.8) (){(a)};
	
		\draw[->, red, ultra thick]  (1.7,-2)--(2.7,-2); 						
	\end{scope}

		\begin{scope}[shift={(6,0)}, rotate=0]

		{\tikzstyle{every node}=[draw ,circle,fill=white, minimum size=0.35cm, inner sep=1pt]						
			\node[draw] at (-0.5,0) (w1){};
			\node[draw] at (0.5,0) (w2){};
			\node[draw] at (-2,-1.5) (u11){$x$};
			\node[draw] at (-2,-2.5) (u12){$y$};
			\node[draw] at (-2+0.5,-2) (u13){$z$};
					
					\node[draw] at (2,-1.5) (u21){};
					\node[draw] at (2,-2.5) (u22){};
					\node[draw] at (2-0.5,-2) (u23){};
					
		\node[draw] at (-2.5,-4) (v41){};
		\node[draw] at (-1.5,-4) (v42){};
		\node[draw] at (-2,-4.5) (v43){};
		\node[draw] at (1.5,-4) (v61){};
		\node[draw] at (2.5,-4) (v62){};
		\node[draw] at (2,-4.5) (v63){};
		\node[draw] at (0-0.5,-4) (v51){};
			\node[draw] at (0+0.5,-4) (v52){};
				\node[draw] at (0,-4+0.5) (v53){};
		
		}
		
		\path[draw,black]
	(w1) edge node[name=la,pos=0.5, above] {\color{blue}  edge in $D$} (w2)				
	(w1) edge node[name=la,pos=0.7, above] {\color{blue} } (u11)		
	(u11) edge node[name=la,pos=0.7, above] {\color{blue} } (u12)		
	(u11) edge node[name=la,pos=0.7, above] {\color{blue} } (u13)		
	(u13) edge node[name=la,pos=0.7, above] {\color{blue} } (u12)		
	(w2) edge node[name=la,pos=0.7, above] {\color{blue} } (u21)		
	(u21) edge node[name=la,pos=0.7, above] {\color{blue} } (u22)		
	(u21) edge node[name=la,pos=0.7, above] {\color{blue} } (u23)		
	(u23) edge node[name=la,pos=0.7, above] {\color{blue} } (u22)		
	(u12) edge node[name=la,pos=0.7, above] {\color{blue} } (v41)		
	(u22) edge node[name=la,pos=0.7, above] {\color{blue} } (v62)		
	(v42) edge node[name=la,pos=0.7, above] {\color{blue} } (v41)		
	(v61) edge node[name=la,pos=0.7, above] {\color{blue} } (v62)		
	(v42) edge node[name=la,pos=0.7, above] {\color{blue} } (v51)		
	(v52) edge node[name=la,pos=0.7, above] {\color{blue} } (v51)		
	(v52) edge node[name=la,pos=0.7, above] {\color{blue} } (v61)		
	(v53) edge node[name=la,pos=0.7, above] {\color{blue} } (v51)
	(v53) edge node[name=la,pos=0.7, above] {\color{blue} } (v52)
	(v43) edge node[name=la,pos=0.7, above] {\color{blue} } (v41)
	(v43) edge node[name=la,pos=0.7, above] {\color{blue} } (v42)
	(v63) edge node[name=la,pos=0.7, above] {\color{blue} } (v61)
	(v63) edge node[name=la,pos=0.7, above] {\color{blue} } (v62)
	; 
	
		\draw[->, red, ultra thick]  (3,-2)--(4, -2); 
		
		\node[] at (0,-4.8) (){(b)};
	\end{scope}

		\begin{scope}[shift={(0,-6)}, rotate=0]

		{\tikzstyle{every node}=[draw ,circle,fill=white, minimum size=0.35cm, inner sep=1pt]						
			\node[draw] at (-0.5,0) (w1){};
			\node[draw] at (0.5,0) (w2){};
			\node[draw] at (-2,-1.5) (u11){};
			\node[draw] at (-2,-2.5) (u12){};
			\node[draw] at (-2+0.5,-2) (u13){$v$};
			
			\node[draw] at (2,-1.5) (u21){};
			\node[draw] at (2,-2.5) (u22){};
			\node[draw] at (2-0.5,-2) (u23){};
			
			\node[draw] at (-2.5,-4) (v41){};
			\node[draw] at (-1.5,-4) (v42){};
			\node[draw] at (1.5,-4) (v61){};
			\node[draw] at (2.5,-4) (v62){};
			\node[draw] at (0-0.5,-4) (v51){};
			\node[draw] at (0+0.5,-4) (v52){};
			\node[draw] at (0,-4+0.5) (v53){};
				\node[draw] at (-2,-4.5) (v43){};
			\node[draw] at (2,-4.5) (v63){};
			
		}

			{\tikzstyle{every node}=[draw ,circle,fill=black, minimum size=0.3cm,
			inner sep=0pt]
			\draw [] (w1) -- (u11) node [midway] (t1){};
			\draw [] (u12) -- (v41) node [midway] (t2){};
			\draw [] (v42) -- (v51) node [midway] (t3){};
			\draw [] (v52) -- (v61) node [midway] (t4){};
			\draw [] (v62) -- (u22) node [midway] (t5){};
			\draw [] (u21) -- (w2) node [midway] (t6){};
			
			\node[draw] at (-2+1,-2) (t7){};
			\node[draw] at (2-1,-2) (t8){};
			\node[draw] at (0,-4+1) (t9){};
			
		}

		\path[draw,black]
		(w1) edge node[name=la,pos=0.5, above] {\color{blue}  edge in $D$} (w2)				
		(w1) edge node[name=la,pos=0.7, above] {\color{blue} } (u11)		
		(u11) edge node[name=la,pos=0.7, above] {\color{blue} } (u12)		
		(u11) edge node[name=la,pos=0.7, above] {\color{blue} } (u13)		
		(u13) edge node[name=la,pos=0.7, above] {\color{blue} } (u12)		
		(w2) edge node[name=la,pos=0.7, above] {\color{blue} } (u21)		
		(u21) edge node[name=la,pos=0.7, above] {\color{blue} } (u22)		
		(u21) edge node[name=la,pos=0.7, above] {\color{blue} } (u23)		
		(u23) edge node[name=la,pos=0.7, above] {\color{blue} } (u22)		
		(u12) edge node[name=la,pos=0.7, above] {\color{blue} } (v41)		
		(u22) edge node[name=la,pos=0.7, above] {\color{blue} } (v62)		
		(v42) edge node[name=la,pos=0.7, above] {\color{blue} } (v41)		
		(v61) edge node[name=la,pos=0.7, above] {\color{blue} } (v62)		
		(v42) edge node[name=la,pos=0.7, above] {\color{blue} } (v51)		
		(v52) edge node[name=la,pos=0.7, above] {\color{blue} } (v51)		
		(v52) edge node[name=la,pos=0.7, above] {\color{blue} } (v61)		
		(v53) edge node[name=la,pos=0.7, above] {\color{blue} } (v51)
		(v53) edge node[name=la,pos=0.7, above] {\color{blue} } (v52)
			(t7) edge node[name=la,pos=0.7, above] {\color{blue} } (u13)
				(t8) edge node[name=la,pos=0.7, above] {\color{blue} } (u23)
					(t9) edge node[name=la,pos=0.7, above] {\color{blue} } (v53)
					(v43) edge node[name=la,pos=0.7, above] {\color{blue} } (v41)
					(v43) edge node[name=la,pos=0.7, above] {\color{blue} } (v42)
					(v63) edge node[name=la,pos=0.7, above] {\color{blue} } (v61)
					(v63) edge node[name=la,pos=0.7, above] {\color{blue} } (v62)
		; 
		
		\draw[->, red, ultra thick]  (3,-2)--(4, -2); 
		
		\node[] at (0,-4.8) (){(c)};
		
		\node[] at (-1,-1.65) (){$v'$};
	\end{scope}

		\begin{scope}[shift={(7,-6)}, rotate=0]

		{\tikzstyle{every node}=[draw ,circle,fill=white, minimum size=0.35cm, inner sep=1pt]						
			\node[draw] at (-0.5,0) (w1){};
			\node[draw] at (0.5,0) (w2){};
			\node[draw] at (-2,-1.5) (u11){};
			\node[draw] at (-2,-2.5) (u12){};
			\node[draw] at (-2+0.5,-2) (u13){};
			
			\node[draw] at (2,-1.5) (u21){};
			\node[draw] at (2,-2.5) (u22){};
			\node[draw] at (2-0.5,-2) (u23){};
			
			\node[draw] at (-2.5,-4) (v41){};
			\node[draw] at (-1.5,-4) (v42){};
			\node[draw] at (1.5,-4) (v61){};
			\node[draw] at (2.5,-4) (v62){};
			\node[draw] at (0-0.5,-4) (v51){};
			\node[draw] at (0+0.5,-4) (v52){};
			\node[draw] at (0,-4+0.5) (v53){};
				\node[draw] at (-2,-4.5) (v43){};
			\node[draw] at (2,-4.5) (v63){};
			
		}

		{\tikzstyle{every node}=[draw ,circle,fill=black, minimum size=0.3cm,
			inner sep=0pt]
			\draw [] (w1) -- (u11) node [midway] (t1){};
			\draw [] (u12) -- (v41) node [midway] (t2){};
			\draw [] (v42) -- (v51) node [midway] (t3){};
			\draw [] (v52) -- (v61) node [midway] (t4){};
			\draw [] (v62) -- (u22) node [midway] (t5){};
			\draw [] (u21) -- (w2) node [midway] (t6){};

		}

			{\tikzstyle{every node}=[draw ,circle,fill=white, minimum size=0.3cm,
			inner sep=0pt]

			\node[draw] at (-2+1,-2) (t7){$x_1$};
			\node[draw] at (2-1,-2) (t8){$x_3$};
			\node[draw] at (0,-4+1) (t9){$x_2$};

		}

			{\tikzstyle{every node}=[draw ,circle,fill=white, minimum size=0.3cm,
			inner sep=0pt]
		
			\node[draw] at (-2+1.5,-2.5) (s1){$s_2$};
			\node[draw] at (2-1.5,-2.5) (s2){$s_3$};
			\node[draw] at (0,-1.3) (s3){$s_1$};
			
		}

		\path[draw,black]
		(w1) edge node[name=la,pos=0.5, above] {\color{blue}  edge in $D$} (w2)				
		(w1) edge node[name=la,pos=0.7, above] {\color{blue} } (u11)		
		(u11) edge node[name=la,pos=0.7, above] {\color{blue} } (u12)		
		(u11) edge node[name=la,pos=0.7, above] {\color{blue} } (u13)		
		(u13) edge node[name=la,pos=0.7, above] {\color{blue} } (u12)		
		(w2) edge node[name=la,pos=0.7, above] {\color{blue} } (u21)		
		(u21) edge node[name=la,pos=0.7, above] {\color{blue} } (u22)		
		(u21) edge node[name=la,pos=0.7, above] {\color{blue} } (u23)		
		(u23) edge node[name=la,pos=0.7, above] {\color{blue} } (u22)		
		(u12) edge node[name=la,pos=0.7, above] {\color{blue} } (v41)		
		(u22) edge node[name=la,pos=0.7, above] {\color{blue} } (v62)		
		(v42) edge node[name=la,pos=0.7, above] {\color{blue} } (v41)		
		(v61) edge node[name=la,pos=0.7, above] {\color{blue} } (v62)		
		(v42) edge node[name=la,pos=0.7, above] {\color{blue} } (v51)		
		(v52) edge node[name=la,pos=0.7, above] {\color{blue} } (v51)		
		(v52) edge node[name=la,pos=0.7, above] {\color{blue} } (v61)		
		(v53) edge node[name=la,pos=0.7, above] {\color{blue} } (v51)
		(v53) edge node[name=la,pos=0.7, above] {\color{blue} } (v52)
		(t7) edge node[name=la,pos=0.7, above] {\color{blue} } (u13)
		(t8) edge node[name=la,pos=0.7, above] {\color{blue} } (u23)
		(t9) edge node[name=la,pos=0.7, above] {\color{blue} } (v53)
		(s3) edge node[name=la,pos=0.7, above] {\color{blue} } (w1)
		(s3) edge node[name=la,pos=0.7, above] {\color{blue} } (w2)
		(s3) edge node[name=la,pos=0.7, above] {\color{blue} } (s2)
		(s3) edge node[name=la,pos=0.7, above] {\color{blue} } (s1)
		(s1) edge node[name=la,pos=0.7, above] {\color{blue} } (s2)
		(t7) edge node[name=la,pos=0.7, above] {\color{blue} } (s3)
		(t7) edge node[name=la,pos=0.7, above] {\color{blue} } (s1)
		(t8) edge node[name=la,pos=0.7, above] {\color{blue} } (s3)
		(t8) edge node[name=la,pos=0.7, above] {\color{blue} } (s2)
		(t9) edge node[name=la,pos=0.7, above] {\color{blue} } (s2)
		(t9) edge node[name=la,pos=0.7, above] {\color{blue} } (s1)
		(s3) edge node[name=la,pos=0.7, above] {\color{blue} } (t1)
		(s3) edge node[name=la,pos=0.7, above] {\color{blue} } (u11)
		(s3) edge node[name=la,pos=0.7, above] {\color{blue} } (u13)
		(s3) edge node[name=la,pos=0.7, above] {\color{blue} } (t6)
		(s3) edge node[name=la,pos=0.7, above] {\color{blue} } (u21)
		(s3) edge node[name=la,pos=0.7, above] {\color{blue} } (u23)
		
		(s1) edge node[name=la,pos=0.7, above] {\color{blue} } (t2)
		(s1) edge node[name=la,pos=0.7, above] {\color{blue} } (u13)
		(s1) edge node[name=la,pos=0.7, above] {\color{blue} } (u12)
		(s1) edge node[name=la,pos=0.7, above] {\color{blue} } (v41)
		(s1) edge node[name=la,pos=0.7, above] {\color{blue} } (v42)
		(s1) edge node[name=la,pos=0.7, above] {\color{blue} } (t3)
		(s1) edge node[name=la,pos=0.7, above] {\color{blue} } (v51)
		(s1) edge node[name=la,pos=0.7, above] {\color{blue} } (v53)

		(s2) edge node[name=la,pos=0.7, above] {\color{blue} } (t4)
		(s2) edge node[name=la,pos=0.7, above] {\color{blue} } (v52)
		(s2) edge node[name=la,pos=0.7, above] {\color{blue} } (v53)
		(s2) edge node[name=la,pos=0.7, above] {\color{blue} } (v61)
		(s2) edge node[name=la,pos=0.7, above] {\color{blue} } (v62)
		(s2) edge node[name=la,pos=0.7, above] {\color{blue} } (t5)
		(s2) edge node[name=la,pos=0.7, above] {\color{blue} } (u22)
		(s2) edge node[name=la,pos=0.7, above] {\color{blue} } (u23)
		
		(v43) edge node[name=la,pos=0.7, above] {\color{blue} } (v41)
		(v43) edge node[name=la,pos=0.7, above] {\color{blue} } (v42)
		(v63) edge node[name=la,pos=0.7, above] {\color{blue} } (v61)
		(v63) edge node[name=la,pos=0.7, above] {\color{blue} } (v62)
		; 
		
		\node[] at (0,-4.8) (){(d)};
		%	\draw[->, red, ultra thick]  (3,-2)--(4, -2); 
	\end{scope}

			\end{tikzpicture}
		\end{center}
		\caption{Constructing $G$ from $G_0$.}
		\label{fig2}
	\end{figure}
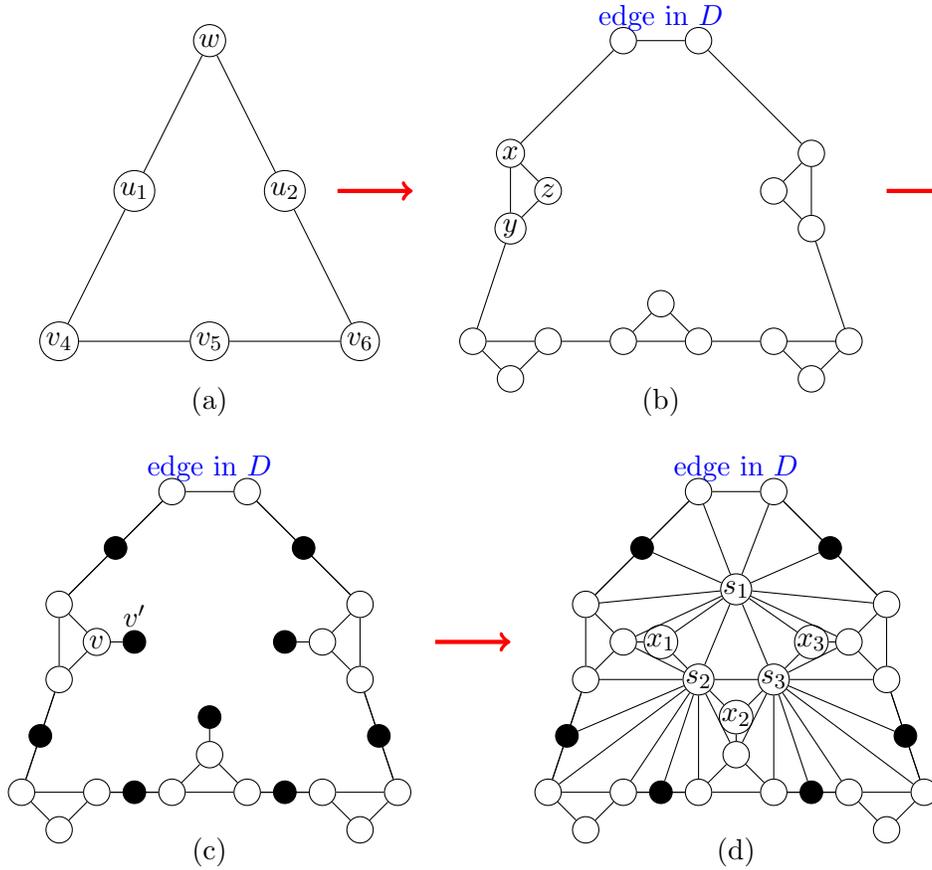

{\noindent {\bf Step 2.2}: In  the graph resulting from Step 2.1, we replace each vertex $v\in V(G_0) \setminus\{w\}$  with a triangle. }  The replacement distinguishes whether $d_{G_0}(v)=2$ or $d_{G_0}(v)=3$, and 
is illustrated below.  
	\begin{enumerate}[(i)]
		\item If $d_{G_0}(v)=3$, then we delete $v$, place a triangle in the position of $v$,  and  add a perfect matching between the three vertices of the triangle and  the three neighbors of $v$ in the current graph; 
		\item If $d_{G_0}(v)=2$, then  $v$ is depicted in black or gray in Figure~\ref{fig1}. We delete $v$, place a triangle $xyzx$ in the position of $v$, 
		add a perfect matching between $\{x,y\}$  and the two neighbors of $v$ in the current graph  such that the vertex $z$ is embedded  inside  the $S$-triangle face  of $G_0$ that is  associated  with $v$.  See Figure~\ref{fig2}(b) for an illustration of the resulting graph of the face boundary $wu_1v_4v_5v_6u_2w$ after this step.   \label{Step-2.2-2}
	\end{enumerate}
	
	The resulting graph from Steps 2.1-2.2  is called $G_1$. 
	 The  faces of $G_1$ incident with a vertex of degree 2 of $G_1$ correspond to 
	the $S$-triangle faces of $G_0$, and are still called $S$-triangle faces of $G_1$.  The triangles used to replace the vertices of $G_0-w$ in this step are called  \emph{$\CC_3$- triangles}.

	\medskip 
	
	{\noindent {\bf Step 3}: The  construction of  $G_2$.} 
	
	\medskip 
	
	We subdivide each edge of $G_1$  joining vertices from two distinct $\CC_3$-triangles,  or  with one from a $\CC_3$-triangle and the other from $D$, and let the set of those new vertices be $T_2$. 
	For each vertex  $v$ of degree 2 in $G_1$, those are the vertices  like $z$ produced in Step 2.2~\eqref{Step-2.2-2}, we add a new vertex $v'$, place $v'$ in the  $S$-triangle face  of $G_1$ that is incident with  $v$, and add the edge $vv'$.   The resulting graph is called $G_2$. 
	Let the set of those new vertices $v'$ be $T_1$. 
	See Figure~\ref{fig2}(c) for an illustration.  The  faces of $G_2$ incident with a vertex of degree 1 of $G_2$ correspond to 
	the $S$-triangle faces of $G_1$, and are still called the $S$-triangle faces of $G_2$. 
	
	Let 
	\begin{equation}
		T=T_1\cup T_2. 
	\end{equation}

	\medskip 

{\noindent {\bf Step 4}: The  construction of  $G$.} 

\medskip 

For each face $F$ of $G_2$, we do the following operations. 
	\begin{enumerate}[(i)]
	\item If $F$ is not an $S$-triangle face, we embed a new vertex inside the the face and join an edge from 
	the new vertex  to all the vertices on the boundary of the face; 
	\item If $F$ is  an $S$-triangle face, let $x_1,x_2, x_3$ be the three vertices of degree 1 incident with the face. We embed a new triangle $s_1s_2s_3s_1$ inside the 
	face,   we first add the edges $x_1s_1, x_1s_2, x_2s_2,x_2s_3, x_3s_3,x_3s_1$, then in the current plane graph, we add edges joining $s_i$ to all the non-adjacent vertices of $s_i$ on the face boundary containing $s_i$, $x_{i-1}$, and $x_i$ for each $i\in [1,3]$, where   $x_0:=x_3$ (triangulate the face).  Denote the set of all new vertices placed in the faces of $G_2$ by $S$. 
	See Figure~\ref{fig2}(d) for an illustration.    
\end{enumerate}

The triangles such as   $s_1s_2s_3s_1$  added in Step 4(ii)  are  called \emph{$S$-triangles}.   We also call $\{s_1,s_2,s_3\}$, the vertex set of an $S$-triangle, an $S$-triangle just for 
notation simplicity. 
If one vertex of an $S$-triangle is adjacent in $G$ to a vertex from $V(D)$, we say that the $S$-triangle is \emph{associated with $D$}. 
Otherwise,  the $S$-triangle is not  associated with $D$. 
The resulting graph from Step 4 is now defined to be $G$.  By the construction, $G$ is a plane triangulation.  Let $p$  and $q$ be  respectively the number of 
vertices of degree 2 and 3 in $G_0$, and let $f_s(G_0)$ be the number of $S$-triangle faces of $G_0$.  
By direct counting and calculations, we have 
\begin{eqnarray*}
p&=&17\times 3+4\times 3=63, \\
q &=& 86-(17+12)=57,\\
n(G_0)&=& p+q+1=121, \\
e(G_0)&=&  \frac{1}{2}(39+2p+3q) =168, \\
f(G_0)&=& 2+e(G_0)-n(G_0)=49, \\
f_s(G_0)&=& 17+4=21,   \\
|S|&=&f(G_0)+2\times f_s(G_0) =91, \\
|T|& =&e(G_0) +3f_s(G_0) =231. 
\end{eqnarray*}

We define some notation before we proceed with the rest proofs.  Let $$U=V(G)\setminus (S\cup T)  \quad \text{and} \quad U_3=U\setminus V(D).$$ 
Vertices from $S$, $T$, and $U$  are called \emph{$S$-vertices}, \emph{$T$-vertices}, and \emph{$U$-vertices}, respectively. 
 
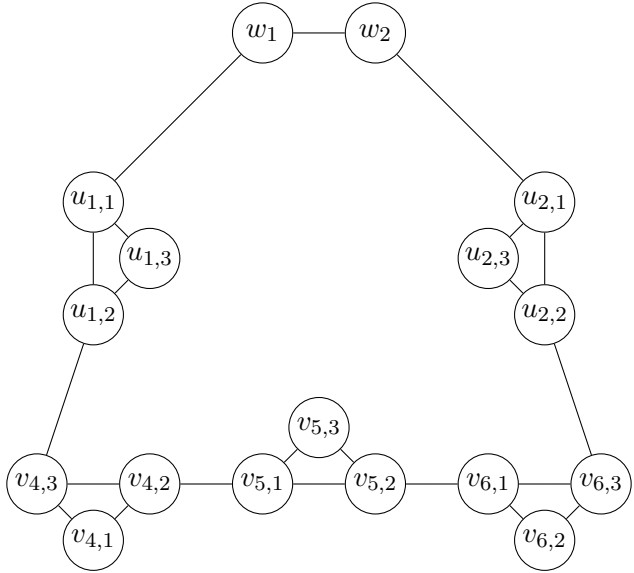
\begin{figure}[!htb]
	
	\begin{center}
		\begin{tikzpicture}[scale=1.5]	
			\usetikzlibrary{calc}
	{\tikzstyle{every node}=[draw ,circle,fill=white, minimum size=0.8cm, inner sep=1pt]						
	\node[draw] at (-0.5,0) (w1){$w_1$};
	\node[draw] at (0.5,0) (w2){$w_2$};
	\node[draw] at (-2,-1.5) (u11){$u_{1,1}$};
	\node[draw] at (-2,-2.5) (u12){$u_{1,2}$};
	\node[draw] at (-2+0.5,-2) (u13){$u_{1,3}$};
	
	\node[draw] at (2,-1.5) (u21){$u_{2,1}$};
	\node[draw] at (2,-2.5) (u22){$u_{2,2}$};
	\node[draw] at (2-0.5,-2) (u23){$u_{2,3}$};
	
	\node[draw] at (-2.5,-4) (v41){$v_{4,3}$};
	\node[draw] at (-1.5,-4) (v42){$v_{4,2}$};
	\node[draw] at (-2,-4.5) (v43){$v_{4,1}$};
	\node[draw] at (1.5,-4) (v61){$v_{6,1}$};
	\node[draw] at (2.5,-4) (v62){$v_{6,3}$};
	\node[draw] at (2,-4.5) (v63){$v_{6,2}$};
	\node[draw] at (0-0.5,-4) (v51){$v_{5,1}$};
	\node[draw] at (0+0.5,-4) (v52){$v_{5,2}$};
	\node[draw] at (0,-4+0.5) (v53){$v_{5,3}$};
	
}

\path[draw,black]
(w1) edge node[name=la,pos=0.5, above] {\color{blue}} (w2)				
(w1) edge node[name=la,pos=0.7, above] {\color{blue} } (u11)		
(u11) edge node[name=la,pos=0.7, above] {\color{blue} } (u12)		
(u11) edge node[name=la,pos=0.7, above] {\color{blue} } (u13)		
(u13) edge node[name=la,pos=0.7, above] {\color{blue} } (u12)		
(w2) edge node[name=la,pos=0.7, above] {\color{blue} } (u21)		
(u21) edge node[name=la,pos=0.7, above] {\color{blue} } (u22)		
(u21) edge node[name=la,pos=0.7, above] {\color{blue} } (u23)		
(u23) edge node[name=la,pos=0.7, above] {\color{blue} } (u22)		
(u12) edge node[name=la,pos=0.7, above] {\color{blue} } (v41)		
(u22) edge node[name=la,pos=0.7, above] {\color{blue} } (v62)		
(v42) edge node[name=la,pos=0.7, above] {\color{blue} } (v41)		
(v61) edge node[name=la,pos=0.7, above] {\color{blue} } (v62)		
(v42) edge node[name=la,pos=0.7, above] {\color{blue} } (v51)		
(v52) edge node[name=la,pos=0.7, above] {\color{blue} } (v51)		
(v52) edge node[name=la,pos=0.7, above] {\color{blue} } (v61)		
(v53) edge node[name=la,pos=0.7, above] {\color{blue} } (v51)
(v53) edge node[name=la,pos=0.7, above] {\color{blue} } (v52)
(v43) edge node[name=la,pos=0.7, above] {\color{blue} } (v41)
(v43) edge node[name=la,pos=0.7, above] {\color{blue} } (v42)
(v63) edge node[name=la,pos=0.7, above] {\color{blue} } (v61)
(v63) edge node[name=la,pos=0.7, above] {\color{blue} } (v62)
; 

%\draw[->, red, ultra thick]  (3,-2)--(4, -2); 
%
%\node[] at (0,-4.8) (){(b)};
\end{tikzpicture}
\end{center}
\caption{The labels of vertices from the $\CC_3$-triangles, where $w_1,w_2\in V(D)$.}
\label{fig2.5}
\end{figure}

For a vertex $v\in V(G_0)\setminus \{w\}$, we let $R(v)$ be the $\CC_3$-triangle that was used to replace $v$ in Step 2.2, and we write $R(v)=v_{i,1}v_{i,2}v_{i,3}v_{i,1}$ if $v=v_i$ for some $i\in [1,86]$,
and $R(v)=u_{i,1}u_{i,2}u_{i,3}u_{i,1}$ if $v=u_i$ for some $i\in [1,34]$.  For each $i\in [1,34]$,  we assume that $u_{i,1}$  and 
a vertex from $D$ have  in $G$ a common neighbor from $T$, and $u_{i,3}$ is embedded inside the  $S$-triangle face of   $G_1$ 
that is corresponding to the $S$-triangle face of $G_0$ associated with  $u_i$.  
%For example, if let $y_{1,1} y_{1,2}y_{1,3}y_{1,1}$ 
%from Figure~\ref{fig4} represent $u_1$, then $y_{1,1} =u_{1,1}, y_{1,2}=u_{1,2}$, and $y_{1,3}=u_{1,3}$. 
Furthermore, we  assume that  for  each vertex $v_{i}$ with $i\in [1,86]$,  if $d_{G_0}(v_i)=3$, then $v_{i,3}$
  and a vertex  $u_{j,2}$ for  some $j\in [1,34] $ or some vertex of $D$
have  in $G$ a common neighbor  from $T$,  and if $d_{G_0}(v_i)=2$, then $v_{i,3}$ is embedded inside the  $S$-triangle face of   $G_1$ 
that is corresponding to the $S$-triangle face of $G_0$ associated with  $v_i$. 
We  let the other two vertices of $R(v_i)$ for $i\in [1,86]$ be labeled such that 
 $$C_1=v_{1,1}v_{1,2} v_{2,1}v_{2,2}  \ldots v_{86,1} v_{86,2} v_{1,1}$$  is a cycle in $G_1$.  For an illustration, see Figure~\ref{fig2.5}. 
The cycle  in $G_2$  obtained from $C_1$ by subdividing each of its edges is denoted by $C$.  
The labels for the vertices  $u_i$'s, vertices from $R(u_i)$, $v_j$'s, and vertices from $R(v_j)$ for $i\in [1,34]$ and $j\in [1,86]$ 
will be fixed throughout the paper.

A path  $P$ in $G$ that  can be denoted as $wt_1u_{i,1}u_{i,2}t_2v_{j,3}$ or $wt_1v_{j,3}$, where $w\in V(D)$,  $t_1,t_2\in T_2$, $i\in [1,34]$, 
and $j\in [1,86]$  is called a  \emph{spoke} of $G$. 
 The former is a \emph{long spoke} while the latter is 
a \emph{short spoke}.  For example, $w_1t_1 x_{1,1}x_{1,2}t_2x_{2,1}$ in Figure~\ref{fig4} is a long spoke.

Let $G^*$ be  obtained from $G_2$ by deleting all the vertices from $T_1$,  deleting $u_{i,3}$ for each $i\in [1,34]$,  and deleting $v_{i,3}$ if $d_{G_0}(v_i)=2$. In other words, by our assumption of the lables of the vertices of $G$, $G^*=G-S-T_1-(N_G(T_1)\cap U)$. 
For a face of $G^*$  with boundary $F$,  we say that $F$ is the \emph{boundary} of the set of $S$-vertices, say  $S^*$,   that are embedded in $G$ inside $F$. 
The subgraph $F\cap C$ is called the \emph{$C$-segment} of $F$. 
 If $S^*$ is associated with $D$, then the boundary of $S^*$ consists of two spokes and one $C$-segment. 
We let $S^*_l$ and $S^*_r$ respectively denote the two sets of $S$-vertices embedded inside the two faces of $G^*$ that share a spoke with $F$ on its left and right, and let $S^*_c$ be the set of $S$-vertices 
embedded inside the face of $G^*$ that share the $C$-segment with  $F$.  
For example, using Figure~\ref{fig4} as an illustration,  if we  let 
$$F=w_1t_1x_{1,1}x_{1,2}t_2x_{2,1}x_{2,2}t_3x_{3,1} x_{3,2}t_4 x_{4,1} x_{4,2} t_5 x_{5,1} x_{5,2} t_6w_2w_1,$$ 
and let  the set of the three $S$-vertices of $G$ embedded inside $F$ be $S^*$,  then the set of the $S$-vertex 
adjacent in $G$ to vertices $ w_2, t_6, x_{5,2}, x_{5,1}, t_5$ is $S_r^*$,   the set of  the $S$-vertex 
adjacent in $G$ to vertices $ w_1, t_1, x_{1,1}, x_{1,2}, t_2$ is $S_l^*$, and the set of the $S$-vertices adjacent to $x_{2,2}, t_3, x_{3,1}, x_{3,2}, t_4$ is $S_c^*$. 
%See Figure~\ref{fig0} for an illustration.  
An $S$-triangle  $S^*$ is internal if 
the boundary of each of $S^*_l$ and $S^*_r$ contains no short spoke. 
%{\color{red}{picture}}

%Two  disjoint subgraphs   $A$ and $B$  are \emph{$T$-adjacent}  in $G$ if there is a vertex from $T$ that is adjacent in $G$  to one vertex from $A$ and one vertex from $B$. 
For an $S$-triangle $s_1s_2s_3s_1$, the three vertices from $T_1$ that each have in $G$ two neighbors from $\{s_1, s_2, s_3\}$ are called the \emph{$T$-vertices associated} with $s_1s_2s_3s_1$. The three $\CC_3$-triangles that each have a vertex adjacent in $G$ to a vertex from the $T$-vertices associated with $s_1s_2s_3s_1$ are called the \emph{$\CC_3$-triangles  associated} with $s_1s_2s_3s_1$.

For a vertex $y\in T$, a neighbor of $y$ from $U$ in $G$ is called a \emph{$U$-neighbor} of $y$ in $G$. 
For a vertex $u\in U$ such that $e_G(u,T) =1$, we let $T(u)$ be the neighbor of $u$ from $T$.

We first show that $G$ has no 2-factor. 
\begin{CLA}\label{claim:(S,T)-barrier}
The pair $(S,T)$ is a barrier of $G$ with  $\delta(S,T)=-2$.  As a consequence, $G$ does not have a 2-factor. 
\end{CLA}

\pf By the construction of $G$,  we know that $T$  is an independent set in $G$, and $G-(S\cup T)$ has no even component.   Thus  $\sum_{y\in T} d_{G-S}(y) =\sum_{k\ge 0}(2k+1) |\CC_{2k+1}|$, and so  
\begin{eqnarray*}
	\delta(S,T) &=& 2|S|-2|T|+ \sum_{k\ge 0}(2k+1) |\CC_{2k+1}|-\sum_{k\ge 0} |\CC_{2k+1}| \\ 
	&=& 2|S|-2|T|+ 2|\CC_3|+38|\CC_{39}| \\
	 &=& 2|S|-2|T|+2(p+q)+38  =182-462+240+38=-2. 
\end{eqnarray*}
Thus by Theorem~\ref{tutte's theorem}, $G$ does not have a 2-factor. 
\qed 

A vertex of  a cutset $W$  of $G$ is said to be \emph{connected} to a component of $G-W$ 
if that vertex is adjacent in $G$ to a vertex  from the component.  
To finish proving  Theorem~\ref{main}, it remains to show that $\tau(G) \ge \frac{3}{2}$.  We will prove this by a  contradiction. 
As $\delta(G)=3$, we have $\tau(G) \le \frac{3}{2}$. Suppose to the contrary that $\tau(G)<\frac{3}{2}$. 
We choose a cutset  $W$ of $G$  such that 
\begin{enumerate}[(a)]
	\item  $h(W):=\frac{3}{2}c(G-W)-|W|$ is as large as possible;  and 
	\item subject to (a), $|W|$ is as large as possible; and 
	\item subject to (b), $|W \cap T|$ is as  small as possible.
\end{enumerate}
Since we assumed $\tau(G)<\frac{3}{2}$, there exists $W\subseteq V(G)$ such that  $h(W) \ge \frac{1}{2}$. 
Also by the constraint  (a) in the choice of $W$,  we have the following fact. 
\begin{FAC}\label{fact:W-vertices-connect-to-components}
	For any $W' \subseteq W$ such that $W\setminus W'$ is a cutset of $G$,  vertices of $W'$ are connected in $G$  to at least $\frac{2}{3}|W'|+1$ 
	components of $G-W$. 
\end{FAC}

%each vertex of $W$ is connected in $G$ to at least two components of $G-W$. 
%We will take a special cutset of $G$, which we call  \emph{irreducible}. 
%A cutset $W$ of $G$ is irreducible if for  any cutset  $W'$  of $G$ with $|W'| \le |W|$, we have $c(G-W')  \le c(G-W)$.  Clearly, if $W$ is a cutset of $G$ such that $\frac{|W|}{c(G-W)}=\tau(G)$, then $W$ is irreducible.  Thus $G$ has an irreducible cutset. 
%We let $W$ be a maximum irreducible cutset of $G$ such that $\frac{|W|}{c(G-W)} <\frac{3}{2}$. 
Our goal is to show that each component of $G-W$ is either a single vertex or   an edge. This structure 
restriction of the components 
forces $W$ to be  consisted of some special vertices. 
Based on that we will show that we actually have   $h(W)<0$, achieving  a contradiction.

\begin{CLA}\label{claim:W-no-intersect-B-C}
	It holds that $W\cap (V(D)\setminus V(A_1)) =\emptyset$.  As a consequence, $D-W$ is connected. 
\end{CLA}

\pf  Recall that $D$ is the replacement graph for the vertex $w$ in $G_0$. 
Suppose to the contrary that $W\cap (V(D)\setminus V(A_1))  \ne \emptyset$. 
Let $W'=W\cap (V(D)\setminus V(A_1)) $. 
As $W'\subseteq V(D)\setminus V(A_1)$ and $G$ 
is a plane graph, we know that for any component $Q$ of $G-W$ with $V(Q)\cap V(D) =\emptyset$,  
we have $E_G(W', Q) =\emptyset$. Hence  the restriction of  the components of $G-W$ to which vertices from $W'$ are connected 
are components of $D-W'$.  Thus each vertex of $W'$  must be 
connected  in $D$ to at least two components of $D-W'$ by Fact~\ref{fact:W-vertices-connect-to-components}.  
Therefore $W'$ is a cutset of $D$, and so 
 $c(D-W') \le \frac{1}{2}|W'|$
as  $D$ is 2-tough by Claim~\ref{fact:D-2-tough}.  
As $h(W)>0$ and  $c(D-W') \le \frac{1}{2}|W'|$,  it follows that 
$|W| \ge |W'|+1$ and $c(G-W) >\frac{2}{3} |W| > \frac{1}{2}|W'|$.  Thus $c(G-W) \ge \lfloor \frac{1}{2}|W'|\rfloor +1$ as $c(G-W)$ is an integer. 
Hence  $c(G-(W\setminus W')) \ge c(G-W) - \lfloor \frac{1}{2}|W'|\rfloor  +2$, 
and so $W\setminus W'$ is a cutset of $G$.  This gives a contradiction to Fact~\ref{fact:W-vertices-connect-to-components}
since  vertices of $W'$ are connected in $G$ to at  most $\frac{1}{2}|W'|$ components of $G-W$. 
The consequence part of the statement is clear as every vertex of $V(A_1)$ has in $D$ a neighbor from $V(A_2)$ by the construction of $D$. 
\qed

\begin{CLA}\label{claim:block-components}
Each component of $G-W$ has no cutvertex. 
\end{CLA}

\pf Suppose to the contrary that a component $Q$ of $G-W$ has a cutvertex, say $x$. 
Then $h(W\cup \{x\})  \ge h(W)+\frac{3}{2} -1>h(W)$, a contradiction to the choice of $W$. 
\qed

\begin{CLA}\label{claim:T-vertex-in-W}
For any $y\in W\cap T$,  we have  $d_G(y) =4$  and $N_G(y)\cap S \subseteq W$.  
\end{CLA}

\pf   Since $y$ is  connected to at least two components of $G-W$, it follows that $G[N_G(y)]$ is not a complete graph. Thus $d_G(y)=4$ 
and so $G[N_G(y)]$ is a 4-cycle.  By the construction of $G$, we have $|N_G(y)\cap S|=|N_G(y)\cap U|=2$. 
Let $N_G(y)=\{u_1,u_2, s_1, s_2\}$ with $u_1, u_2\in U$ and $s_1, s_2\in S$. Then by Fact~\ref{fact:W-vertices-connect-to-components}, we must have  that either $u_1,u_2\in W$
and  $s_1$ and $s_2$ are  separated  in two different components of $G-W$, or $s_1,s_2\in W$
and  $u_1$ and $u_2$ are  separated  in two different components of $G-W$. Suppose  to the contrary that $u_1,u_2\in W$
and  $s_1$ and $s_2$ are  separated  in two different components of $G-W$ is the case.  For $i\in [1,2]$, let $Q_i$ 
be the odd component of $G-(S\cup T)$ containing $u_i$. By the construction of $G$,      each of $s_1$ and $s_2$ has 
two neighbors in $G$ from $V(Q_i)$,  and $s_1$ and $s_2$ have in $G$ exactly one common neighbor from $V(Q_i)$. 
 If $Q_i$ is a triangle, then $Q_i-W$ has at most one component; if $Q_i=D$, then again $Q_i-W$ is connected by Claim~\ref{claim:W-no-intersect-B-C}. 
As $s_1$ and $s_2$ are  separated  in two different components of $G-W$, 
we must have $|W\cap V(Q_i)\cap N_G(\{s_1,s_2\})| \ge 2$ for each $i\in [1,2]$.   
Let $|(W\cap N_G(\{s_1,s_2\})\cap V(Q_1)) \cup (W\cap N_G(\{s_1,s_2\})\cap V(Q_2))|  =j$, where $j\in [4,6]$.   
As  $|N_G(s_1)\cap V(Q_i)|=|N_G(s_2)\cap V(Q_i)|=2$ and 
$s_1$ and $s_2$ have in $G$ exactly one common neighbor from $Q_i$,  
it follows that 
 all vertices from  $W':=\{y\} \cup (W\cap N_G(\{s_1,s_2\})\cap V(Q_1)) \cup (W\cap N_G(\{s_1,s_2\})\cap V(Q_2))$ 
are connected in $G$ to at most $2+(j-4)=j-2$ components of $G-W$: deleting  two vertices from $W\cap N_G(\{s_1,s_2\})\cap V(Q_1)$ 
and two vertices from $W\cap N_G(\{s_1,s_2\})\cap V(Q_2)$ can create at most two components in $G-W$, and each third deletion of a vertex from $W\cap N_G(\{s_1,s_2\})\cap V(Q_1)$  
or $W\cap N_G(\{s_1,s_2\})\cap V(Q_2)$ can create at most one more component either by the adjacencies of vertices of a $\CC_3$-triangle or by Claim~\ref{claim:W-no-intersect-B-C} if one of $Q_i$ is $D$.  As $h(W)>0$ and $|W| \ge j+1$, we have $c(G-W) > \frac{2}{3}(j+1)$. Since $\frac{2}{3}(j+1) -(j-2) =\frac{8}{3} -\frac{j}{3} >0$, we know 
that $W\setminus W'$ is a cutset of $G$. 
This shows a contradiction to Fact~\ref{fact:W-vertices-connect-to-components} as $j-2<\frac{2}{3}(j+1)+1$ when $j\in [4,6]$. 
\qed

\begin{CLA}\label{claim:T-vertices-adj-A}
	For any $y\in T\cap W$,  $y$ is connected  in $G$ to exactly two components of $G-W$, and both of the components are trivial. 
\end{CLA}

\pf  If $y\in W$, then we have  $d_G(y)=4$ and $N_G(y)\cap S \subseteq W$  by Claim~\ref{claim:T-vertex-in-W}. Thus 
$y$ is connected  in $G$ to exactly two components of $G-W$ by Fact~\ref{fact:W-vertices-connect-to-components}. 
Furthermore,  the two components 
of $G-W$ connected to  $y$ in $G$ respectively contain  the two neighbors of $y$ from $U$. 
Let $N_G(y)\cap U=\{u_1,u_2\}$ and assume that the component of $G-W$ that contains $u_1$ is not trivial. 
Then  $(W\setminus \{y\}) \cup \{u_1\}$ 
is a cutset of $G$ with the same size as $W$ and $h((W\setminus \{y\}) \cup \{u_1\}) =h(W)$. However $(W\setminus \{y\}) \cup \{u_1\}$ 
contains less vertices from $T$ than $W$ does, a contradiction to the choice of $W$. 
\qed

\begin{CLA}\label{claim:C-vertex-in-W}
	For any $y\in W\cap U$,  we have  $N_G(y)\cap S \subseteq W$.  
\end{CLA}

\pf  Let $Q$ be the component of $G-(S\cup T)$ containing $y$.  
By the construction of $G$ and Claim~\ref{claim:W-no-intersect-B-C}, we know that $y$ 
is adjacent in $G$ to a vertex from $T$, two vertices from $S$, and some vertices from  $Q$. 
As $Q-W$ is connected, Fact~\ref{fact:W-vertices-connect-to-components} implies that $y$ is connected in $G$ to exactly two components of $G-W$, where either the two components contain vertices of $Q-W$ and $N_G(y)\cap T$ respectively, or  the two components respectively  contain one distinct vertex from $N_G(y)\cap S$. Let us suppose instead that 
$N_G(y)\cap S \not\subseteq W$. Then   we must have $N_G(y)\cap T \subseteq W$. 
Let $z$ be the vertex from $N_G(y)\cap T$. As $N_G(z)\cap S= N_G(y)\cap S$, Claim~\ref{claim:T-vertex-in-W} implies  $N_G(y)\cap S \subseteq W$, a contradiction. 
\qed 

\begin{CLA}\label{claim:outerface-component}
	Let $Q$ be a 2-connected component of $G-W$, and $F$ be the boundary of a  face of $Q$ such that $E_G(F, W) \ne \emptyset$. Then  $V(F)\cap S=\emptyset$. 
%	\begin{enumerate}[(i)]
%		\item If $V(F)\cap S \ne \emptyset$, then every vertex of $V(F)\cap S$ is from an $S$-triangle not associated with $D$. 
%	\end{enumerate}
\end{CLA}

\pf  Suppose to the contrary that $V(F)\cap S \ne \emptyset$. 
Let $s_1\in V(F)\cap S$. As  $s_1\in V(F)$, $E_G(V(F), W) \ne \emptyset$, and  $G$ is a plane triangulation, we    have  $E_G(s_1, W) \ne \emptyset$.  Then by Claims~\ref{claim:T-vertex-in-W} and~\ref{claim:C-vertex-in-W}, we have $N_G(s_1) \cap (T\cup U) \subseteq Q$. 
As a consequence, $|V(Q)| \ge |N_G(s_1) \cap (T\cup U)| \ge 6$.  Since $N_G(s_1) \cap (T\cup U) \subseteq Q$ and $E_G(s_1, W) \ne \emptyset$, 
it follows that  $s_1$ is one vertex from an 
$S$-triangle, say $s_1s_2s_3s_1$.   Let $t_1, t_2, t_3$ 
be the three vertices from $T$ such that $t_is_i, t_is_{i+1}\in E(G)$ for each $i\in [1,3]$, where $s_4:=s_1$.  Let $u_1, u_2, u_3\in U$ such that $t_iu_i \in E(G)$ for $i\in [1,3]$. 

If $s_2, s_3\in W$, then we have $t_1, t_3\in V(Q)$.   
Then $c(Q-\{s_1, u_1, u_3\} )\ge 3$. Therefore, $W \cup\{s_1, u_1, u_3\}$ is a cutset of $G-W$ with  $h(W \cup\{s_1, u_1, u_3\}) \ge h(W)$, a contradiction to the choice of $W$. 
Thus we assume that $|W\cap \{s_2, s_3\}| \le 1$.  Since $E_G(s_1, W)=E_G(s_1, W\cap \{s_2,s_3\}) \ne \emptyset$, 
 we may assume, without loss of generality, that $s_2\in V(Q)$ and $s_3 \in W$.  Then  $(N_G(s_1) \cup N_G(s_2))\setminus \{s_3\} \subseteq V(Q)$ by Claims~\ref{claim:T-vertex-in-W} and~\ref{claim:C-vertex-in-W}. 
This further implies that $t_1, t_2, t_3\in V(Q)$ by Claim~\ref{claim:T-vertex-in-W}. 
Let $F'$ be the boundary of $s_1s_2s_3s_1$, where recall that $F'$  is defined as a face boundary of $G^*$ before Claim~\ref{claim:(S,T)-barrier}.  Then we know that one component of $F'-W$ intersects 
with $Q$, and the rest are paths.   

%By Claim~\ref{claim:block-components},  it must be the case that the other vertices of $S$ incident with $F'$ are all contained in 
%$W$.  If all of those $S$-vertices is contained in $W$, then $F'-W$ is connected, a contradiction to Fact~\ref{fact:W-vertices-connect-to-components}. 
%Let $s'$ be the $S$-vertex  connected to a path component of $F'-W$ in $G$,  and suppose  that in the next face of $G_1$ in the clockwise direction, the $S$-vertex 
%adjacent to that segment of $F'$ is not contained in $W$. Then  the corresponding component of $G-W$ of the segment of $F'$ has to be either a single vertex or an edge, 
%or 2-connected. 

Then we must have $|V(F')\cap W\cap N_G(s_3)| \ge 5$ and  $c(F'-(W\cap N_G(s_3))) \ge 5$. For otherwise,  let $|V(F')\cap W|  =j$ for  $j\in [2,4]$.  As $F'$ is a cycle and so is 1-tough, we have $c(F'-W)  \le j$. 
This implies that vertices of $V(F')\cap W\cap N_G(s_3)$ are connected in $G$ to at most $j$ components of $G-W$. 
However  $j<\frac{2}{3} (j+1)+1$, contradicting Fact~\ref{fact:W-vertices-connect-to-components}.  Thus $|V(F')\cap W\cap N_G(s_3)| \ge 5$ and  $c(F'-(W\cap N_G(s_3))) \ge 5$. 

\begin{figure}[!htb]

	\begin{center}
		\begin{tikzpicture}[scale=1]	
			\usetikzlibrary{calc}
			\begin{scope}[shift={(0,0)}, rotate=0]		
%drawing vertices
\node [draw, circle, label={[shift={(0.2,0.2)}] left:$v_{9,1}$}] at  (-4.5,2.5) (v91) {};
\node [draw, circle, label={[shift={(-0.2,0.2)}]right:$v_{9,2}$}] at (-3.5,2.5) (v92) {};
\node [draw, circle, label={above,align=right:$v_{9,3}$}] at (-4, 3) (v93) {};			
\node [draw, circle, label={[shift={(0.2,0.2)}]left:$v_{10,1}$}] at (-1.5,2.5) (v101) {};
\node [draw, circle, label={[shift={(-0.2,0.2)}]right:$v_{10,2}$}] at (-0.5,2.5) (v102) {};
\node [draw, circle, label={above,align=right:$v_{10,3}$}] at  (-1, 3) (v103) {};				
\node [draw, circle, label={[shift={(0.2,0.2)}]left:$v_{11,1}$}] at (-1.5+3,2.5) (v111) {};
\node [draw, circle, label={[shift={(-0.2,0.2)}]right:$v_{11,2}$}] at (-0.5+3,2.5) (v112) {};
\node [draw, circle, label={[shift={(-0.15,0)}]right,align=right:$v_{11,3}$}] at (-1+3, 2) (v113) {};					
\node [draw, circle, label={[shift={(0.2,0.2)}]left:$v_{12,1}$}] at (-1.5+6,2.5) (v121) {};
\node [draw, circle, label={[shift={(-0.2,0.2)}]right:$v_{12,2}$}] at (-0.5+6,2.5) (v122) {};
\node [draw, circle, label={above,align=right:$v_{12,3}$}] at  (-1+6, 3) (v123) {};		
\node [draw, circle, label={[shift={(0.2,0.2)}]left:$v_{13,1}$}] at  (-1.5+9,2.5) (v131) {};
\node [draw, circle, label={[shift={(-0.2,0.2)}]right:$v_{13,2}$}] at (-0.5+9,2.5) (v132) {};
\node [draw, circle, label={above,align=right:$v_{13,3}$}] at (-1+9, 3) (v133) {};	

\node [draw, circle, label={[shift={(0.15,-0.2)}]left,align=right:$v_{8,1}$}] at (-4.5, -.5) (v81) {};	
\node [draw, circle, label={[shift={(0.15,0.2)}]left,align=right:$v_{8,2}$}] at (-4.5, 0.5) (v82) {};	
\node [draw, circle, label={left,align=right:$v_{8,3}$}] at (-5, 0) (v83) {};	

\node [draw, circle, label={[shift={(0.15,-0.2)}]left,align=right:$v_{7,1}$}] at (-4.5, -3.5) (v71) {};	
\node [draw, circle, label={[shift={(0.15,0.2)}]left,align=right:$v_{7,2}$}] at (-4.5, -2.5) (v72) {};	
\node [draw, circle, label={left,align=right:$v_{7,3}$}] at (-5, -3) (v73) {};

\node [draw, circle, label={[shift={(0.15,-0.2)}]left,align=right:$v_{20,2}$}] at (-4.5, -5.5) (v202) {};	
\node [draw, circle, label={[shift={(0,0)}]above,align=right:$v_{20,1}$}] at (-3.5, -5.5) (v201) {};	
\node [draw, circle, label={below,align=right:$v_{20,3}$}] at (-4, -6) (v203) {};	

\node [draw, circle, label={[shift={(0,0)}]below,align=right:$v_{19,2}$}] at (-4.5+3, -5.5) (v192) {};	
\node [draw, circle, label={[shift={(0,0)}]below,align=right:$v_{19,1}$}] at (-3.5+3, -5.5) (v191) {};	
\node [draw, circle, label={[shift={(-0.28, -0.05)}]above,align=right:$v_{19,3}$}] at (-4+3, -5) (v193) {};

\node [draw, circle, label={[shift={(0,0)}]above,align=right:$v_{18,2}$}] at (-4.5+6, -5.5) (v182) {};	
\node [draw, circle, label={[shift={(0,0)}]above,align=right:$v_{18,1}$}] at (-3.5+6, -5.5) (v181) {};	
\node [draw, circle, label={below,align=right:$v_{18,3}$}] at (-4+6, -6) (v183) {};

\node [draw, circle, label={[shift={(0,0)}]above,align=right:$v_{17,2}$}] at (-4.5+9, -5.5) (v172) {};	
\node [draw, circle, label={[shift={(0,0)}]above,align=right:$v_{17,1}$}] at (-3.5+9, -5.5) (v171) {};	
\node [draw, circle, label={below,align=right:$v_{17,3}$}] at (-4+9, -6) (v173) {};

\node [draw, circle, label={[shift={(0,0)}]above,align=right:$v_{16,2}$}] at (-4.5+12, -5.5) (v162) {};	
\node [draw, circle, label={[shift={(-0.15,-0.2)}]right,align=right:$v_{16,1}$}] at (-3.5+12, -5.5) (v161) {};	
\node [draw, circle, label={below,align=right:$v_{16,3}$}] at (-4+12, -6) (v163) {};

\node [draw, circle, label={[shift={(-0.15,-0.2)}]right,align=right:$v_{14,2}$}] at (-4.5+13, -.5) (v142) {};	
\node [draw, circle, label={[shift={(-0.15,0.2)}]right,align=right:$v_{14,1}$}] at (-4.5+13, 0.5) (v141) {};	
\node [draw, circle, label={right,align=right:$v_{14,3}$}] at (-4+13, 0) (v143) {};	

\node [draw, circle, label={[shift={(-0.15,-0.2)}]right,align=right:$v_{15,2}$}] at (-4.5+13, -3.5) (v152) {};	
\node [draw, circle, label={[shift={(-0.15,0.2)}]right,align=right:$v_{15,1}$}] at (-4.5+13, -2.5) (v151) {};	
\node [draw, circle, label={[shift={(0.15,-0.2)}]left,align=right:$v_{15,3}$}] at (-5+13, -3) (v153) {};

\node [draw, circle, label={right,align=right:$t_a$}] at (-1+3, 0) (ta) {};	
\node [draw, circle, label={left,align=right:$t_b$}] at (-0.5, -3) (tb) {};	
\node [draw, circle, label={below,align=right:$t_c$}] at (4.5, -3) (tc) {};	

\node [draw, circle, label={left,align=right:$a$}] at (-0.5, -1.5) (a) {};	
\node [draw, circle, label={right,align=right:$b$}] at (2, -4.5) (b) {};	
\node [draw, circle, label={right,align=right:$c$}] at (4.5, -1.5) (c) {};	

\path[draw,black]
(v91) edge node[name=la,pos=0.5, above] {\color{blue}  } (v92)				
(v91) edge node[name=la,pos=0.7, above] {\color{blue} } (v93)		
(v92) edge node[name=la,pos=0.7, above] {\color{blue} } (v93)		
(v101) edge node[name=la,pos=0.7, above] {\color{blue} } (v102)		
(v101) edge node[name=la,pos=0.7, above] {\color{blue} } (v103)		
(v102) edge node[name=la,pos=0.7, above] {\color{blue} } (v103)		
(v111) edge node[name=la,pos=0.7, above] {\color{blue} } (v112)		
(v111) edge node[name=la,pos=0.7, above] {\color{blue} } (v113)		
(v112) edge node[name=la,pos=0.7, above] {\color{blue} } (v113)		
(v121) edge node[name=la,pos=0.7, above] {\color{blue} } (v122)	
(v121) edge node[name=la,pos=0.7, above] {\color{blue} } (v123)	
(v122) edge node[name=la,pos=0.7, above] {\color{blue} } (v123)	
(v131) edge node[name=la,pos=0.7, above] {\color{blue} } (v132)	
(v131) edge node[name=la,pos=0.7, above] {\color{blue} } (v133)	
(v132) edge node[name=la,pos=0.7, above] {\color{blue} } (v133)	
(v92) edge node[name=la,pos=0.7, above] {\color{blue} } (v101)		
(v102) edge node[name=la,pos=0.7, above] {\color{blue} } (v111)		
(v112) edge node[name=la,pos=0.7, above] {\color{blue} } (v121)		
(v122) edge node[name=la,pos=0.7, above] {\color{blue} } (v131)		
(v91) edge node[name=la,pos=0.7, above] {\color{blue} } (v82)		
(v82) edge node[name=la,pos=0.7, above] {\color{blue} } (v81)		
(v82) edge node[name=la,pos=0.7, above] {\color{blue} } (v83)		
(v81) edge node[name=la,pos=0.7, above] {\color{blue} } (v83)		
(v81) edge node[name=la,pos=0.7, above] {\color{blue} } (v72)		
(v72) edge node[name=la,pos=0.7, above] {\color{blue} } (v71)		
(v72) edge node[name=la,pos=0.7, above] {\color{blue} } (v73)		
(v71) edge node[name=la,pos=0.7, above] {\color{blue} } (v73)		
(v71) edge node[name=la,pos=0.7, above] {\color{blue} } (v202)		
(v202) edge node[name=la,pos=0.7, above] {\color{blue} } (v201)		
(v202) edge node[name=la,pos=0.7, above] {\color{blue} } (v203)		
(v201) edge node[name=la,pos=0.7, above] {\color{blue} } (v203)		

(v201) edge node[name=la,pos=0.7, above] {\color{blue} } (v192)		
(v192) edge node[name=la,pos=0.7, above] {\color{blue} } (v191)		
(v192) edge node[name=la,pos=0.7, above] {\color{blue} } (v193)		
(v191) edge node[name=la,pos=0.7, above] {\color{blue} } (v193)		

(v191) edge node[name=la,pos=0.7, above] {\color{blue} } (v182)		
(v182) edge node[name=la,pos=0.7, above] {\color{blue} } (v181)		
(v182) edge node[name=la,pos=0.7, above] {\color{blue} } (v183)		
(v181) edge node[name=la,pos=0.7, above] {\color{blue} } (v183)		

(v181) edge node[name=la,pos=0.7, above] {\color{blue} } (v172)		
(v172) edge node[name=la,pos=0.7, above] {\color{blue} } (v171)		
(v172) edge node[name=la,pos=0.7, above] {\color{blue} } (v173)		
(v171) edge node[name=la,pos=0.7, above] {\color{blue} } (v173)		

(v171) edge node[name=la,pos=0.7, above] {\color{blue} } (v162)		
(v162) edge node[name=la,pos=0.7, above] {\color{blue} } (v161)		
(v162) edge node[name=la,pos=0.7, above] {\color{blue} } (v163)		
(v161) edge node[name=la,pos=0.7, above] {\color{blue} } (v163)		

(v161) edge node[name=la,pos=0.7, above] {\color{blue} } (v152)		
(v152) edge node[name=la,pos=0.7, above] {\color{blue} } (v151)		
(v152) edge node[name=la,pos=0.7, above] {\color{blue} } (v153)		
(v151) edge node[name=la,pos=0.7, above] {\color{blue} } (v153)		

(v151) edge node[name=la,pos=0.7, above] {\color{blue} } (v142)		
(v142) edge node[name=la,pos=0.7, above] {\color{blue} } (v141)		
(v142) edge node[name=la,pos=0.7, above] {\color{blue} } (v143)		
(v141) edge node[name=la,pos=0.7, above] {\color{blue} } (v143)		

(v141) edge node[name=la,pos=0.7, above] {\color{blue} } (v132)		
(v132) edge node[name=la,pos=0.7, above] {\color{blue} } (v131)		
(v132) edge node[name=la,pos=0.7, above] {\color{blue} } (v133)		
(v131) edge node[name=la,pos=0.7, above] {\color{blue} } (v133)		

(v113) edge node[name=la,pos=0.7, above] {\color{blue} } (ta)		
(v193) edge node[name=la,pos=0.7, above] {\color{blue} } (tb)		
(v153) edge node[name=la,pos=0.7, above] {\color{blue} } (tc)		

(a) edge node[name=la,pos=0.7, above] {\color{blue} } (b)		
(a) edge node[name=la,pos=0.7, above] {\color{blue} } (c)		
(b) edge node[name=la,pos=0.7, above] {\color{blue} } (c)		

(ta) edge node[name=la,pos=0.7, above] {\color{blue} } (a)		
(ta) edge node[name=la,pos=0.7, above] {\color{blue} } (c)		
(tb) edge node[name=la,pos=0.7, above] {\color{blue} } (a)		
(tb) edge node[name=la,pos=0.7, above] {\color{blue} } (b)		
(tc) edge node[name=la,pos=0.7, above] {\color{blue} } (b)		
(tc) edge node[name=la,pos=0.7, above] {\color{blue} } (c)		
;

	{\tikzstyle{every node}=[draw ,circle,fill=white, minimum size=0.3cm,
	inner sep=0pt]
	\draw [] (v202) -- (v71) node [midway,label={[shift={(0.05, 0)}]left,align=right:$t_1$}] (t1){};
	
	\draw [] (v72) -- (v81) node [midway,label={[shift={(0.05, 0)}]left,align=right:$t_2$}] (t2){};
	\draw [] (v82) -- (v91) node [midway,label={[shift={(0.05, 0)}]left,align=right:$t_3$}] (t3){};
		\draw [] (v92) -- (v101) node [midway,label={[shift={(0.05, 0)}]below,align=right:$t_4$}] (t4){};
	\draw [] (v102) -- (v111) node [midway,label={[shift={(0.05, 0)}]below,align=right:$t_5$}] (t5){};
		\draw [] (v112) -- (v121) node [midway,label={[shift={(0.05, 0)}]below,align=right:$t_6$}] (t6){};
		\draw [] (v122) -- (v131) node [midway,label={[shift={(0.05, 0)}]below,align=right:$t_7$}] (t7){};
		\draw [] (v132) -- (v141) node [midway,label={[shift={(0.05, 0)}]right,align=right:$t_8$}] (t8){};
		\draw [] (v142) -- (v151) node [midway,label={[shift={(0.05, 0)}]right,align=right:$t_9$}] (t9){};
		\draw [] (v152) -- (v161) node [midway,label={[shift={(0.05, 0)}]right,align=right:$t_{10}$}] (t10){};
		
		\draw [] (v162) -- (v171) node [midway,label={[shift={(0.05, 0)}]above,align=right:$t_{11}$}] (t11){};
		\draw [] (v172) -- (v181) node [midway,label={[shift={(0.05, 0)}]above,align=right:$t_{12}$}] (t12){};
		\draw [] (v182) -- (v191) node [midway,label={[shift={(0.05, 0)}]above,align=right:$t_{13}$}] (t13){};
		\draw [] (v192) -- (v201) node [midway,label={[shift={(0.05, 0)}]above,align=right:$t_{14}$}] (t14){};

}

%\path[draw,black]
%
%(a) edge node[name=la,pos=0.7, above] {\color{blue} } (v193)		
%(a) edge node[name=la,pos=0.7, above] {\color{blue} } (v192)		
%(a) edge node[name=la,pos=0.7, above] {\color{blue} } (t14)		
%(a) edge node[name=la,pos=0.7, above] {\color{blue} } (v201)		
%(a) edge node[name=la,pos=0.7, above] {\color{blue} } (v202)		
%(a) edge node[name=la,pos=0.7, above] {\color{blue} } (t1)		
%;

\end{scope}

\end{tikzpicture}
\end{center}
\caption{An $S$-triangle $abca$ and  the neighbors of the vertices $a,b,c$. The edges joining  $a$, $b$  and $c$ to other vertices for triangulating the three faces of length more than 3 are omitted in the drawing.}
\label{fig3}
\end{figure}
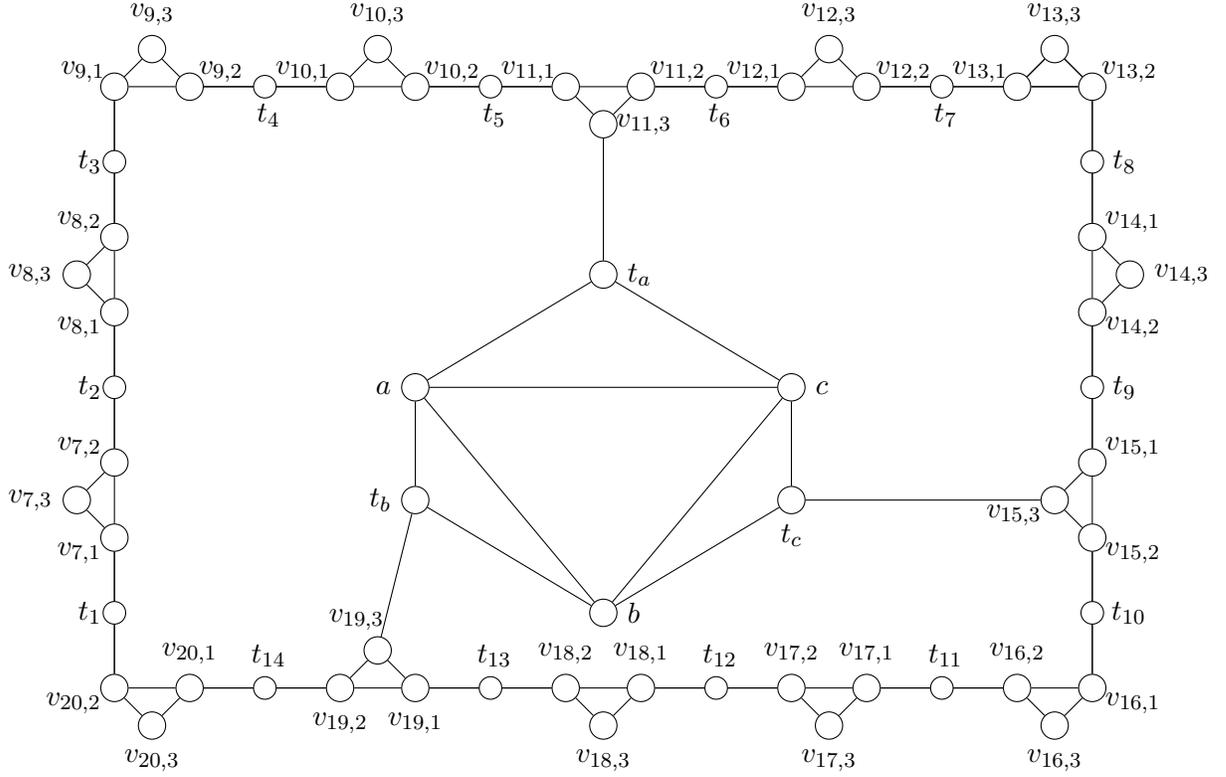

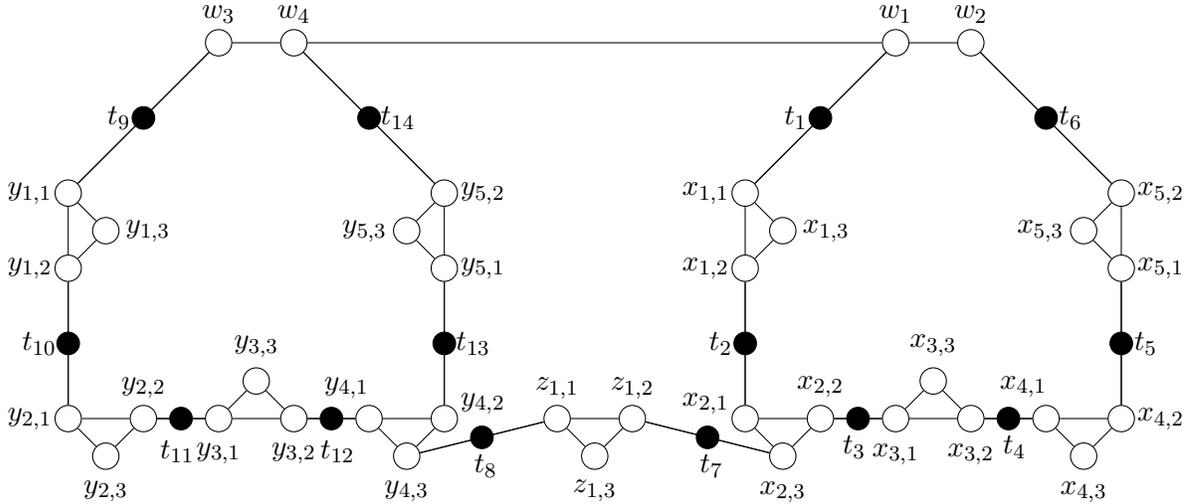
\begin{figure}[!htb]
	\begin{center}
		\begin{tikzpicture}[scale=1]	
			\usetikzlibrary{calc}

				{\tikzstyle{every node}=[draw ,circle,fill=white, minimum size=0.35cm, inner sep=1pt]
					\node [draw, circle, label={[shift={(0,-0.1)}] above:$w_3$}] at  (-0.5,2) (w3) {};						
					\node [draw, circle, label={[shift={(0,-0.1)}] above:$w_4$}] at  (0.5,2) (w4) {};		
					\node [draw, circle, label={[shift={(0.05,0)}] left:$y_{1,1}$}] at  (-2.5,0) (y11) {};		
					\node [draw, circle, label={[shift={(0.05,0)}] left:$y_{1,2}$}] at  (-2.5,-1) (y12) {};		
					\node [draw, circle, label={[shift={(0,0)}] right:$y_{1,3}$}] at  (-2,-0.5) (y13) {};		
					
					\node [draw, circle, label={[shift={(-0.05,0)}] right:$y_{5,2}$}] at  (2.5,0) (y22) {};		
					\node [draw, circle, label={[shift={(-0.05,0)}] right:$y_{5,1}$}] at  (2.5,-1) (y21) {};		
					\node [draw, circle, label={[shift={(0,0)}] left:$y_{5,3}$}] at  (2,-0.5) (y23) {};

					\node [draw, circle, label={[shift={(0.05,0)}] left:$y_{2,1}$}] at  (-2.5,-3) (w41) {};		
					\node [draw, circle, label={[shift={(0,-0.15)}] above:$y_{2,2}$}] at  (-1.5,-3) (w42) {};		
					\node [draw, circle, label={[shift={(0,0.1)}] below:$y_{2,3}$}] at  (-2,-3.5) (w43) {};		
					
					\node [draw, circle, label={[shift={(0,0.1)}] below:$y_{3,1}$}] at  (-2.5+2,-3) (w51) {};		
					\node [draw, circle, label={[shift={(0,0.1)}] below:$y_{3,2}$}] at  (-1.5+2,-3) (w52) {};		
					\node [draw, circle, label={[shift={(0,-0.1)}] above:$y_{3,3}$}] at  (-2+2,-2.5) (w53) {};

					\node [draw, circle, label={[shift={(-0.3,-0.15)}]above :$y_{4,1}$}] at  (-2.5+4,-3) (w61) {};		
					\node [draw, circle, label={[shift={(-0.05,0.2)}] right:$y_{4,2}$}] at  (-1.5+4,-3) (w62) {};		
					\node [draw, circle, label={[shift={(0,0.1)}] below:$y_{4,3}$}] at  (-2+4,-3.5) (w63) {};

						\node [draw, circle, label={[shift={(0,-0.15)}]above :$z_{1,1}$}] at  (4,-3) (z11) {};		
					\node [draw, circle, label={[shift={(0,-0.15)}] above:$z_{1,2}$}] at  (5,-3) (z12) {};		
					\node [draw, circle, label={[shift={(0,0.1)}] below:$z_{1,3}$}] at  (4.5,-3.5) (z13) {};		
					
				}
				
				\path[draw,black]
				(w3) edge node[name=la,pos=0.5, above] {\color{blue} } (w4)	
				(w3) edge node[name=la,pos=0.5, above] {\color{blue} } (y11)			
				(y11) edge node[name=la,pos=0.5, above] {\color{blue} } (y12)
				(y12) edge node[name=la,pos=0.5, above] {\color{blue} } (y13)
				(y11) edge node[name=la,pos=0.5, above] {\color{blue} } (y13)
				(y12) edge node[name=la,pos=0.5, above] {\color{blue} } (w41)
				(w42) edge node[name=la,pos=0.5, above] {\color{blue} } (w41)
				(w42) edge node[name=la,pos=0.5, above] {\color{blue} } (w43)
				(w41) edge node[name=la,pos=0.5, above] {\color{blue} } (w43)
				(w42) edge node[name=la,pos=0.5, above] {\color{blue} } (w51)
				(w52) edge node[name=la,pos=0.5, above] {\color{blue} } (w51)
				(w53) edge node[name=la,pos=0.5, above] {\color{blue} } (w51)
				(w52) edge node[name=la,pos=0.5, above] {\color{blue} } (w53)
				(w52) edge node[name=la,pos=0.5, above] {\color{blue} } (w61)
				(w62) edge node[name=la,pos=0.5, above] {\color{blue} } (w61)
				(w63) edge node[name=la,pos=0.5, above] {\color{blue} } (w61)
				(w63) edge node[name=la,pos=0.5, above] {\color{blue} } (w62)
				(w62) edge node[name=la,pos=0.5, above] {\color{blue} } (y21)
				(y21) edge node[name=la,pos=0.5, above] {\color{blue} } (y22)
				(y21) edge node[name=la,pos=0.5, above] {\color{blue} } (y23)
				(y22) edge node[name=la,pos=0.5, above] {\color{blue} } (y23)
				(y22) edge node[name=la,pos=0.5, above] {\color{blue} } (w4)
				; 
				
				{\tikzstyle{every node}=[draw ,circle,fill=black, minimum size=0.3cm,
					inner sep=0pt]
					\draw [] (w3) -- (y11) node [midway,label={[shift={(0.05, 0)}]left,align=right:$t_9$}] (t9){};
					
					\draw [] (y12) -- (w41) node [midway,label={[shift={(0.05, 0)}]left,align=right:$t_{10}$}] (t10){};
					\draw [] (w42) -- (w51) node [midway,label={[shift={(-0.05, 0)}]below,align=right:$t_{11}$}] (t11){};
					\draw [] (w52) -- (w61) node [midway,label={[shift={(0.08, 0)}]below,align=right:$t_{12}$}] (t12){};
					\draw [] (w62) -- (y21) node [midway,label={[shift={(-0.05, 0)}]right,align=right:$t_{13}$}] (t13){};
					\draw [] (y22) -- (w4) node [midway,label={[shift={(-0.05, 0)}]right,align=right:$t_{14}$}] (t14){};
					
				}

			{\tikzstyle{every node}=[draw ,circle,fill=white, minimum size=0.35cm, inner sep=1pt]
				\node [draw, circle, label={[shift={(0,-0.1)}] above:$w_1$}] at  (-0.5+9,2) (w1) {};						
				\node [draw, circle, label={[shift={(0,-0.1)}] above:$w_2$}] at  (0.5+9,2) (w2) {};		
				\node [draw, circle, label={[shift={(0.05,0)}] left:$x_{1,1}$}] at  (-2.5+9,0) (u11) {};		
				\node [draw, circle, label={[shift={(0.05,0)}] left:$x_{1,2}$}] at  (-2.5+9,-1) (u12) {};		
				\node [draw, circle, label={[shift={(0,0)}] right:$x_{1,3}$}] at  (-2+9,-0.5) (u13) {};		
				
				\node [draw, circle, label={[shift={(-0.05,0)}] right:$x_{5,2}$}] at  (2.5+9,0) (u22) {};		
				\node [draw, circle, label={[shift={(-0.05,0)}] right:$x_{5,1}$}] at  (2.5+9,-1) (u21) {};		
				\node [draw, circle, label={[shift={(0,0)}] left:$x_{5,3}$}] at  (2+9,-0.5) (u23) {};

				\node [draw, circle, label={[shift={(0.05,0.2)}] left:$x_{2,1}$}] at  (-2.5+9,-3) (v41) {};		
				\node [draw, circle, label={[shift={(0,-0.15)}] above:$x_{2,2}$}] at  (-1.5+9,-3) (v42) {};		
				\node [draw, circle, label={[shift={(0,0.1)}] below:$x_{2,3}$}] at  (-2+9,-3.5) (v43) {};		
				
				\node [draw, circle, label={[shift={(0,0.1)}] below:$x_{3,1}$}] at  (-2.5+2+9,-3) (v51) {};		
				\node [draw, circle, label={[shift={(0,0.1)}] below:$x_{3,2}$}] at  (-1.5+2+9,-3) (v52) {};		
				\node [draw, circle, label={[shift={(0,-0.1)}] above:$x_{3,3}$}] at  (-2+2+9,-2.5) (v53) {};

				\node [draw, circle, label={[shift={(-0.3,-0.15)}]above :$x_{4,1}$}] at  (-2.5+4+9,-3) (v61) {};		
				\node [draw, circle, label={[shift={(-0.05,0)}] right:$x_{4,2}$}] at  (-1.5+4+9,-3) (v62) {};		
				\node [draw, circle, label={[shift={(0,0.1)}] below:$x_{4,3}$}] at  (-2+4+9,-3.5) (v63) {};

			}
			
			\path[draw,black]
			(w1) edge node[name=la,pos=0.5, above] {\color{blue} } (w2)	
			(w1) edge node[name=la,pos=0.5, above] {\color{blue} } (u11)			
			(u11) edge node[name=la,pos=0.5, above] {\color{blue} } (u12)
			(u12) edge node[name=la,pos=0.5, above] {\color{blue} } (u13)
			(u11) edge node[name=la,pos=0.5, above] {\color{blue} } (u13)
			(u12) edge node[name=la,pos=0.5, above] {\color{blue} } (v41)
			(v42) edge node[name=la,pos=0.5, above] {\color{blue} } (v41)
			(v42) edge node[name=la,pos=0.5, above] {\color{blue} } (v43)
			(v41) edge node[name=la,pos=0.5, above] {\color{blue} } (v43)
			(v42) edge node[name=la,pos=0.5, above] {\color{blue} } (v51)
			(v52) edge node[name=la,pos=0.5, above] {\color{blue} } (v51)
			(v53) edge node[name=la,pos=0.5, above] {\color{blue} } (v51)
			(v52) edge node[name=la,pos=0.5, above] {\color{blue} } (v53)
			(v52) edge node[name=la,pos=0.5, above] {\color{blue} } (v61)
			(v62) edge node[name=la,pos=0.5, above] {\color{blue} } (v61)
			(v63) edge node[name=la,pos=0.5, above] {\color{blue} } (v61)
			(v63) edge node[name=la,pos=0.5, above] {\color{blue} } (v62)
			(v62) edge node[name=la,pos=0.5, above] {\color{blue} } (u21)
			(u21) edge node[name=la,pos=0.5, above] {\color{blue} } (u22)
			(u21) edge node[name=la,pos=0.5, above] {\color{blue} } (u23)
			(u22) edge node[name=la,pos=0.5, above] {\color{blue} } (u23)
			(u22) edge node[name=la,pos=0.5, above] {\color{blue} } (w2)
			; 
			
			{\tikzstyle{every node}=[draw ,circle,fill=black, minimum size=0.3cm,
				inner sep=0pt]
				\draw [] (w1) -- (u11) node [midway,label={[shift={(0.05, 0)}]left,align=right:$t_1$}] (t1){};
				
				\draw [] (u12) -- (v41) node [midway,label={[shift={(0.05, 0)}]left,align=right:$t_2$}] (t2){};
				\draw [] (v42) -- (v51) node [midway,label={[shift={(-0.05, 0)}]below,align=right:$t_3$}] (t3){};
				\draw [] (v52) -- (v61) node [midway,label={[shift={(0.08, 0)}]below,align=right:$t_4$}] (t4){};
				\draw [] (v62) -- (u21) node [midway,label={[shift={(-0.05, 0)}]right,align=right:$t_5$}] (t5){};
				\draw [] (u22) -- (w2) node [midway,label={[shift={(-0.05, 0)}]right,align=right:$t_6$}] (t6){};
				
			}						
			
		%	\node[] at (0,-4.8+9) (){(a)};
			
				\path[draw,black]
				(w4) edge node[name=la,pos=0.5, above] {\color{blue} } (w1)	
			(z11) edge node[name=la,pos=0.5, above] {\color{blue} } (w63)	
			(z12) edge node[name=la,pos=0.5, above] {\color{blue} } (v43)
			(z12) edge node[name=la,pos=0.5, above] {\color{blue} } (z11)
			(z12) edge node[name=la,pos=0.5, above] {\color{blue} } (z13)
			(z11) edge node[name=la,pos=0.5, above] {\color{blue} } (z13);

		{\tikzstyle{every node}=[draw ,circle,fill=black, minimum size=0.3cm,
			inner sep=0pt]
			\draw [] (z12) -- (v43) node [midway,label={[shift={(0.05, 0)}]below,align=right:$t_7$}] (t7){};
			
			\draw [] (z11) -- (w63) node [midway,label={[shift={(0.05, 0)}]below,align=right:$t_8$}] (t8){};
			
		}

		\end{tikzpicture}
	\end{center}
	\caption{A subgraph of  $G$ containing four spokes, where each black vertex $t_i$ with $i\in [1,14]$ is a $T$-vertex and $w_1,w_2, w_3, w_4\in V(D)$. }
	\label{fig4}
\end{figure}

When  $s_1s_2s_3s_1$  is associated with $D$, by the construction of $G$, we know that $s_3$ is adjacent   in $G$ to  at most  6 vertices  from $V(F')$ that are not contained in 
$Q$. Thus it is impossible to   get $c(F'-(W\cap N_G(s_3))) \ge 5$ as $F'$ is 1-tough.  Therefore $s_1s_2s_3s_1$  is not associated with $D$. 
Since for any two $S$-triangles $S_1$ and $S_2$ that are not associated with $D$, we have $G[N_G(S_1)] \cong G[N_G(S_2)]$, 
 we 
may assume that the component graph  $F_0'$ of $F'$ is $v_7v_8\ldots v_{20}v_7$.  We label vertices on $F'$ as shown in Figure~\ref{fig3} for easily referring to them.  
For each $i\in [7,20]$, the triangle $v_{i,1}v_{i,2}v_{i,3}v_{i,1}$ is the replace graph of $v_i$ from $F_0'$,   $\{a,b,c\}=\{s_1,s_2, s_3\} $, and $\{t_a, t_b, t_c, t_1,  t_2, \ldots, t_{14}\} \subseteq T$. 
Since $s_1,s_2\in V(Q)$ and  $G[N_G(b)] \cong G[N_G(c)]$, we may assume that $s_1=c$.  

Consider first that $s_2=a$ and $s_3=b$.  Let the  $S$-vertex  that is adjacent  in $G$ to each vertex from $\{v_{20,3}, v_{20,1}, t_{14}, v_{19,2}, v_{19,1}, t_{13}, v_{18,2}, v_{18,3}\}$ be $s^*$.  
%and let  the  $S$-vertex  that is adjacent in $G$ to each vertex from $\{v_{14,3}, v_{14,2}, t_{9}, v_{15,1}, v_{15,2}, t_{10}, v_{16,2}, v_{16,3}\}$ be $s'$. 
If $W\cap \{v_{18,2}, v_{19,1}\} \ne \emptyset$, then we must have $s^*\in W$ by Claim~\ref{claim:C-vertex-in-W}.  Then  we have  $c(Q-\{a,c, v_{11,3}, v_{19,3}, v_{15,3}, t_{14}\})\ge 5$.  
However, $|W\cup \{a,c, v_{11,3}, v_{19,3}, v_{15,3}, t_{14}\}|>|W|$ and 
$h(W\cup \{a,c, v_{11,3}, v_{19,3}, v_{15,3}, t_{14}\}) \ge h(W)$, a contradiction to the choice of $W$.  Therefore we assume that $v_{18,2}, v_{19,1} \in V(Q)$
 and so $t_{13} \in V(Q)$ by Fact~\ref{fact:W-vertices-connect-to-components}. 
 Similarly, we have $v_{15,2}, v_{16,1}, t_{10} \in V(Q)$. Thus the other components of $F'-(W\cap N_G(s_3))$ not containing a vertex from $V(Q)$ will possibly only contain vertices  from  $\{t_{12}, v_{17,2}, v_{17,1}, t_{11}\}$. 
 However, the maximum number of components we can have  
 in $F'[\{t_{12}, v_{17,2}, v_{17,1}, t_{11}\}]$  by deleting its vertices is 2. Thus $c(F'-(W\cap N_G(s_3))) \le 3$, 
 a contradiction  to the assumption that  $c(F'-(W\cap N_G(s_3))) \ge 5$. 

%Then we must have $|W\cap \{s^*,s'\}| \ge  1$. For otherwise,  we have $(N_G(s^*) \cup N_G(s'))\cap (T\cup U) \subseteq V(Q)$ by Claims~\ref{claim:T-vertex-in-W} and~\ref{claim:C-vertex-in-W}.  
%Thus, 
%we  may assume that $s^*\in W$ (the argument for  $s'\in W$ follows from  the same logic). 
%If $v_{19,1}\not\in W$, then we have $t_{13} \not\in W$ by Claim~\ref{claim:T-vertices-adj-A}. 
%As $Q$ is 2-connected, then either the   $S$-vertex,  say $s''$,    that is adjacent  in $G$ to vertices from $\{v_{17,2}, v_{17,3}, v_{18,1},  v_{18,3}, t_{12}\}$  must be contained in $Q$.  Then $N_G(s'')\cap (T\cup U) \subseteq V(Q)$ by Claims~\ref{claim:T-vertex-in-W} and~\ref{claim:C-vertex-in-W}. 
%Then again we get  $c(F'-(W\cap N_G(s_3))) \le 3$, 
%a contradiction to  the assumption that  $c(F'-(W\cap N_G(s_3))) \ge 5$.
%Therefore  $v_{19,1}\in W$. Then  we have  $c(Q-\{a,c, v_{11,3}, v_{19,3}, v_{15,3}, t_{14}\})\ge 5$.  
%However, $|W\cup \{a,c, v_{11,3}, v_{19,3}, v_{15,3}, t_{14}\}|>|W|$ and 
%$h(W\cup \{a,c, v_{11,3}, v_{19,3}, v_{15,3}, t_{14}\}) \ge h(W)$, a contradiction to the choice of $W$. 
 
 Next we assume that $s_2=b$ and $s_3=a$. 
 Suppose first that $W\cap \{v_{19,2}, v_{20,1}\} \ne \emptyset$. Then the  $S$-vertex 
 adjacent in $G$ to both $v_{19,2}$ and $v_{20,1}$ is contained in $W$ by Claim~\ref{claim:C-vertex-in-W}. 
 Then  we have  $c(Q-\{b,c, v_{11,3}, v_{19,3}, v_{15,3}, t_{13}\})\ge 5$. 
However, $|W\cup \{b,c, v_{11,3}, v_{19,3}, v_{15,3}, t_{13}\}|>|W|$ and 
$h(W\cup \{b,c, v_{11,3}, v_{19,3}, v_{15,3}, t_{13}\}) \ge h(W)$, a contradiction to the choice of $W$.  Thus $W\cap \{v_{19,2}, v_{20,1}\} = \emptyset$. 
Thus $v_{19,2} \in V(Q)$ and so $t_{14} \in V(Q)$ by Claim~\ref{claim:T-vertices-adj-A}. 
Hence $v_{20,1} \in V(Q)$ as $W\cap \{v_{19,2}, v_{20,1}\} = \emptyset$. 
 Suppose then  that $W\cap \{v_{20,2}, v_{7,1}\} \ne \emptyset$. Then the $S$-vertex 
 adjacent in $G$  to  both $v_{20,2}$ and $v_{7,1}$ is contained in $W$ by Claim~\ref{claim:C-vertex-in-W}. 
 Let $s^*$ be the  $S$-vertex  adjacent to both $t_{13}$ and $t_{14}$. 
 Then  have  $c(Q-\{b,c, v_{11,3}, v_{19,3}, v_{15,3}, t_{14}, v_{19,1}, v_{18,2}, s^*\})\ge 7$. 
 However, $|W\cup \{b,c, v_{11,3}, v_{19,3}, v_{15,3}, t_{14}, v_{19,1}, v_{18,2}, s^*\}>|W|$ and 
 $$h(W\cup \{b,c, v_{11,3}, v_{19,3}, v_{15,3}, t_{14}, v_{19,1}, v_{18,2}, s^*\}) \ge h(W),$$ a contradiction to the choice of $W$.  Thus $W\cap \{v_{20,2}, v_{7,1}\} = \emptyset$. 
 Thus $v_{20,2} \in V(Q)$ and so $t_{1} \in V(Q)$ by Claim~\ref{claim:T-vertices-adj-A}. 
 Hence $v_{7,1} \in V(Q)$ as $W\cap \{v_{20,2}, v_{7,1}\} = \emptyset$.   
 Symmetrically, we also must have   $v_{11,1}, t_5, v_{10,2}, v_{10,1}, t_4, v_{92} \in V(Q)$. 
 Thus the other components of $F'-(W\cap N_G(s_3))$ not containing a vertex from $V(Q)$ will possibly only contain vertices  from  $\{t_2, v_{8,1}, v_{8,2}, t_3\}$. 
 However, the maximum number of components we can have  
 in $F'[\{t_2, v_{8,1}, v_{8,2}, t_3\}]$  by deleting its vertices is 2. Thus $c(F'-(W\cap N_G(s_3))) \le 3$, 
 a contradiction  to the assumption that  $c(F'-(W\cap N_G(s_3))) \ge 5$. 
 
Therefore  $V(F)\cap S=\emptyset$, proving the statement.  
\qed

%In the component graph $G_0$, 
%we call a path connecting  the vertex $w$ with a vertex of $C_0$ \emph{a spoke} of $G_0$. 
%The corresponding path  in $G$ is called a spoke of $G$. 
\begin{CLA}\label{claim:spoke}
Let $P$ be a spoke of $G$ with endvertices as $x$ and $y$, where $x\in V(D)$ and $y \in V(C)$.  If it holds that $x\not\in W$ or that  $y\not\in W$  but  $y$ is not a component of $G-W$, 
then $P$ is contained in a component of $G-W$. 
\end{CLA}

\pf Suppose first that $x\not\in W$.  Then the vertex, say $t_x$, from $T$ adjacent to $x$ is not contained in $W$ by Claim~\ref{claim:T-vertices-adj-A}.  
Let $s_1,s_2\in S$ be the two neighbors of $t_x$ in $G$. If  $\{s_1, s_2\} \not\subseteq  W$, say $s_1\not\in W$, then  none  of  the $S$-vertices with the same boundary as  $s_1$
is contained in $W$ by Claim~\ref{claim:outerface-component}.  As a consequence, 
 we have $V(P)\cap W =\emptyset$ by Claims~\ref{claim:T-vertex-in-W} and~\ref{claim:C-vertex-in-W}.  Thus we assume $s_1,s_2\in W$.
As a consequence,  all the  $S$-vertices with the same boundary as  $s_1$ or $s_2$ 
are  contained in $W$ by Claim~\ref{claim:outerface-component}.  
  Then as every component of $G-W$ has no cutvertex by Claim~\ref{claim:block-components},   we know 
  that $P$ is entirely contained in a component of $G-W$. 
The argument for  $y\not\in W$ follows from the same idea.  
\qed

%\begin{CLA}\label{claim:spoke-and-C-end}
%	Let $P$ be a spoke of $G$ with $x$ as  its endvertex from $C$. Suppose $P$ is contained in a component, say $Q$, of $G-W$. Then $Q$ contains $x$ and also a vertex  from a $\CC_3$-triangle that is 
%	$T$-adjacent to $x$ and is intersecting $C$. 
%\end{CLA}
%
%\pf Let $S^*$ be the set of $S$-vertices incident with $P$ in one of the two faces, say $F$,   of $G_2$  such that  $|V(F) \cap V(P)| \ge 2$. If $S^* \not\subseteq W$ then $S^*\cap W=\emptyset$ by Claim~\ref{claim:outerface-component}, then the $C$-section of $F$ is included in $Q$ by 
%Claims~\ref{claim:T-vertex-in-W} and~\ref{claim:C-vertex-in-W}.  If both sets of  the $S$-vertices inside the two faces of $G_2$ that have at least two vertices in common with  $P$ are contained in $W$,  then at least one of the $\CC_3$-triangles that   intersects $C$ and is $T$-adjacent to $x$  has a vertex in common with $Q$, 
% as $Q$ is 2-connected by Claim~\ref{claim:block-components}. 
%\qed 

%\begin{CLA}\label{claim:T-neighbor-of-U-vertex}
%Let $abca$ be a replacement triangle such that some of its vertices is contained in a 2-connected component  $Q$ of $G-W$, and it has a vertex a vertex, say $a$, on a face boundary $F$ of 
%$Q$ such that $e_G(F,W) \ge 1$. Then $a\in W$, $T(a)\in W$, and $b,c\in Q$. 
%\end{CLA}
%
%\pf Suppose some vertices of $W$ are embedded in the exterior of $F$. 
%Then all the $S$ vertices embedded in the exterior of $F$ are contained in $W$. 
%As $
%
%\qed 

\begin{CLA}\label{claim:one-2-connected-component}
It holds that $G-W$ has at most one 2-connected component,  which is the component  containing  a vertex  of $D$. 
\end{CLA}

\pf 
 Otherwise, let $Q$ be a 2-connected component of $G-W$ that contains no vertex of $D$, and  let $Q_0$ be the component graph of $Q$. 
 Then $Q_0$ is also 2-connected by Claim~\ref{claim:outerface-component}:  the boundary  $F$ of each face  of $Q$ for which $E_G(F,W) \ne \emptyset$   is a cycle containing no vertex of $S$ and so the component graph $F_0$ of $F$ is a cycle in $G_0$.    Since $Q$ does not contain any vertex of $D$, it follows that  any graph from $\CC$ that intersects $Q$ is a $\CC_3$-triangle. 
  Let $F$ be a face of $Q$ with  a vertex $u\in V(F)\cap U$   satisfying  $e_G(u,W\cap U) \ge 1$. Let  $u'\in W\cap U$  be a vertex for which 
 $e_G(u, u') =1$, and let   $t_{u'} \in T$ such that $u't_{u'} \in E(G)$. 
 Then as $u' \in W$, we have $N_G(u')\cap S \subseteq W$ by Claim~\ref{claim:C-vertex-in-W} and $t_{u'}$ being a 
 component of $G-W$ by Fact~\ref{fact:W-vertices-connect-to-components}. Let $u^*$ be the other vertex for which $u,u'$ together form a $\CC_3$-triangle. Then we must have $u^*\in V(Q)$. 
 For otherwise, we have $u^*\in W$ and so the $S$-vertex that is adjacent to both $u^*$ and $u$ in $G$ is also contained in $W$ by Claim~\ref{claim:C-vertex-in-W}. 
 This implies that the only neighbor of $u$ in $Q$ is from $N_G(u) \cap T$ as $N_G(u')\cap S \subseteq W$.  Since  $|N_G(u)\cap T| =1$, we get a contradiction to $Q$ being 2-connected.  
 Thus $u, u^*\in V(Q)$ and  thus  $N_G(\{u,u^*\})\cap T \subseteq V(Q)$ by Claim~\ref{claim:T-vertices-adj-A}.  Hence if a $\CC$-triangle intersects both $W$ and $Q$, then the vertex  of $Q_0$ that is corresponding to 
 this $\CC$-triangle has degree 2 in $Q_0$. 
 
% Then we must have $|V(Q)\cap \{u,u'u^*\}| \le 2$. For otherwise, we have $N_G(\{u,u'u^*\})\cap T \subseteq V(Q)$ by Claim~\ref{claim:T-vertices-adj-A}
% 
% By Claim~\ref{claim:spoke},  
%we know that  vertices from $V(C)\cap V(Q)$
%are consecutive on $C$.  
Let $b$ be the number of $\CC_3$-triangles that intersects both $W$ and $Q$,   and  $f_s$ be the number of $S$-triangles in $Q$.   We show that $Q_0$ has exactly $3f_s+b$ vertices of degree 2.  Let $v\in V(Q_0)$ such that $d_{Q_0}(v)=2$, and let $R$ be the $\CC_3$-triangle
corresponding to $v$. If $V(R)\cap W \ne \emptyset$ and $V(R)\cap V(Q) \ne \emptyset$, then 
$v$ corresponds to a $\CC_3$-triangle that intersects both $W$ and $Q$ and so there are at most $b$ such degree 2 vertices in $Q$. 
Thus we assume $R \subseteq Q$. Then $v$ must also be a vertex of $G_0$ that is of degree 2 in $G_0$. By Claim~\ref{claim:T-vertices-adj-A}, all the three vertices, say $x,y,z$,  from $T$ that are adjacent in $G$ to vertices from $R$ are contained in $Q$ as well. Since $Q$ is 2-connected, some vertices from 
the $S$-triangle  $S^*$ for which some of  its vertices are adjacent in $G$ to  $x, y, z$ are contained in $Q$ too.  By Claim~\ref{claim:outerface-component}, the entire $S$-triangle  $S^*$ is contained in $Q$.   Let $F$ be the  boundary of the face of $G_2$ such that $S^*$ is embedded inside $F$. 
By Claims~\ref{claim:T-vertex-in-W} and~\ref{claim:C-vertex-in-W},  we know that  $F\subseteq Q$. As $F$ contains vertices from all the three $\CC_3$-triangles  associated with $S^*$ where one of them is $R$, it follows that the three  $\CC_3$-triangles in $Q$ that are   associated with  $S^*$ are all  contained in $Q$.  Thus $v$ is corresponding to one $\CC_3$-triangle that is associated with $S^*$ and so there are at most $3f_s$ such degree 2 vertices in $Q$. 
The above argument shows  that the number of vertices of degree 2 in $Q_0$ is at most $b+3f_s$. 

On the other hand, by the argument from the first paragraph of this proof, we know that 
every $\CC_3$-triangle $R$ that intersects both $W$ and $Q$ corresponds to a vertex of degree 2 in $Q_0$. Furthermore,  for any $S$-triangle  $S^*$ with $S^*\subseteq V(Q)$, let $F$ be the  boundary of the face of $G_2$ such that $S^*$ is embedded inside $F$.   Then we have 
$F\subseteq Q$ by Claims~\ref{claim:T-vertex-in-W} and~\ref{claim:C-vertex-in-W}.  
As $F$ contains all the three $\CC_3$-triangles  associated with $S^*$ together with  the neighbors of the vertices of $S^*$ from $T$ in $G$, it follows that $S^*$ corresponds to three vertices of degree 2 of $Q_0$. Thus $Q_0$ has at least $b+3f_s$ vertices of degree 2.  Combining this with the  assertion that the number of vertices of degree 2 in $Q_0$ is at most $b+3f_s$, we know that $Q_0$ has exactly $3f_s+b$ vertices of degree 2. 
Thus we have  $e(Q_0)  = \frac{1}{2}(3n(Q_0)-3f_s-b)$. 
Consequently by Euler's formula, $f(Q_0)=e(Q_0)-n(Q_0)+2 = 0.5n(Q_0)-1.5f_s-0.5b+2$. 

Let $S^*$ be the set of  $S$-vertices that are embedded inside a face of $G_2$ with boundary $F$.   By Claim~\ref{claim:outerface-component}, we have either $S^*\subseteq V(Q)$ or $S^*\cap V(Q) =\emptyset$. If $S^*\subseteq V(Q)$, then we have $F\subseteq Q$ by Claims~\ref{claim:T-vertex-in-W} and~\ref{claim:C-vertex-in-W}.  Thus $Q_0$ has a face 
whose boundary is the component graph of $F$.  Therefore we have  $|V(Q)\cap S| \le f(Q_0)-1+2f_s$ as at least one face of $Q$ whose boundary has vertices adjacent in $G$ to vertices from $W$ and so 
there is no   $S$-vertex  embedded in $Q$  inside that face.

We will contract a cutset $W_Q$ of $Q$ for which $h(W\cup W_Q) \ge h(W)$ to get a contradiction to the choice of $W$. 
Let $M_0$ be a maximum matching in the component graph $Q_0$ of $Q$.  For each $uv\in M_0$,  there exists $x_u\in R(u)$ and $x_v\in R(v)$
such that $x_u$ and $x_v$ are  both adjacent in $G$ 
to a vertex  from $T$, where recall that $R(u)$ is the $\CC_3$-triangle corresponding to $u$.   
We call $x_u$  a \emph{representative vertex}
of $R(u)$.  As $M_0$ is a matching, each  component  from 
$\CC_3$   either has no representative vertex or has 
a unique representative vertex. 

Let $R\in \CC_{3}$ such that $V(R)\cap V(Q) \ne \emptyset$.  By the argument in the first paragraph in the proof of Claim~\ref{claim:one-2-connected-component}, we have that either $R\subseteq Q$ or 
$|V(R) \cap V(Q)| =2$. 
If $R$ has a representative vertex, say  $x$,  
let $W_R \subseteq (V(R) \cap V(Q))\setminus \{x\}$  be the  set of 
two   vertices (if  $R\subseteq Q$) or one  vertex (if $|V(R) \cap V(Q)|=2$) such that $e_G(R-W_R, T)=e_G(x,T)$. 
Otherwise,  let $W_R \subseteq V(R)$  be a  set of 
two   vertices (if  $R\subseteq Q$) or one  vertex (if $|V(R) \cap V(Q)|=2$) such that  $e_G(R-W_R, T)=1$. 

Note that for any two representative vertices $x_1$ and $x_2$,  $x_1$
and $x_2$ are adjacent in $G$ to the same vertex 
from $T$. We let $T^*$ be the set of all 
these vertices from $T$ that are adjacent in $G$
to a representative vertex of  a  $\CC_3$-triangle  that intersects $Q$.

%, and so $f(Q_0)=e(Q_0)-n(Q_0)+2 =0.5n(Q_0)-1.5f_0-\frac{1}{2}b+2$.  

%For any $D\in \CC_{2k+1}^1$ for some $k\ge 2$, by Claim~\ref{claim:large-odd-comp},  we let $S_D \subseteq V(D)$  be a  set of  $\lfloor \frac{4k+2}{3} \rfloor$ vertices
% such that $D-S_D$ is set of isolated vertices. 

%For any $D\in \CC_{2k+1}^2$ for some $k\ge 2$, we let $S_D \subseteq V(D)$  be the  set of 
%$2k+1$ vertices such that $e_G(D-S_D, T)=0$. 
%Note that $D-S_D$ is a graph with at least one vertex. 

%As there are $2m$ representative vertices, we have $|T^*| = m$. 
%Then we have $|M_0| \ge \frac{1}{2}n(Q_0)$. 
Let $W_Q$ be the set that consists of all vertices in $S\cap V(Q)$,  $T^*$, and $W_R$ for all $R\in \CC_3$ such that $|V(R)\cap V(Q)| \ge 2$.  
Thus $|W_Q| \le f(Q_0)-1+2f_s+2n(Q_0)-b+|T^*|$. 
Each component of $Q-W_Q$ is either a vertex from $T$,  or an edge consisting of a vertex from $T$ and a vertex from $U$, or a vertex from $U$. For the last 
case, the vertex from $U$ corresponds to an endvertex of an edge of $M_0$. Since  the three vertices from $T_1\cap V(Q)$ that are adjacent in $G$ to vertices from every $S$-triangle of $Q$ are contained in $Q$ by Claim~\ref{claim:T-vertex-in-W}, we get 
 $c(Q-W_Q) =e(Q_0)-|T^*|+ 3f_s +2|T^*|$.

 Since $Q$ contains no vertex of $D$ and so contains no spoke of $G$,  it follows that every $\CC_3$-triangle $R$ such that $V(R)\cap V(Q) \ne  \emptyset$ and that $R$ contains a vertex, say $x$,  of a spoke  of $G$  satisfies 
 the property that $x\in W$ by Claim~\ref{claim:spoke}.   
 As a consequence, all the $S$-triangles  that are associated with $D$ and are  adjacent in $G$ to a vertex from $V(Q)\cap V(C)$ are all contained in $W$ 
 by Claims~\ref{claim:C-vertex-in-W} and~\ref{claim:outerface-component}.  Thus $Q_0$  contains  only vertices of $C_0$.   
If a vertex  $t\in T_2\cap V(C)$ is not contained in $Q$, then the two $S$-vertices, say $s_1$ and $s_2$,  adjacent to $t$ in $G$ must be contained in $W$  by Claims~\ref{claim:T-vertex-in-W} (if $t\in W$) or~\ref{claim:T-vertices-adj-A} (if $t\not\in W$).  Thus the two set of $S$-vertices with the same boundary  as $s_1$ or $s_2$ are all contained in $W$ by Claim~\ref{claim:outerface-component}. 
As $Q_0$ is 2-connected,  it then follows that $Q_0\cap C_0$ is connected. This, together with the fact that $Q_0$  contains  only vertices of $C_0$ gives 
 $|T^*| = \lfloor \frac{1}{2} n(Q_0) \rfloor $. 
% Since $Q_0$  contains  only vertices of $C_0$ and every $\CC_3$-triangle $R$ such that $V(R)\cap V(Q) \ne  \emptyset$ and that $R$ contains a vertex, say $x$,  of a spoke  of $G$  satisfies 
% the property that $x\in W$ by Claim~\ref{claim:spoke}, it follows that $b\ge \frac{1}{3}n(Q_0)$ by the construction of $Q_0$. 
  Letting $\alpha=0$   if $n(G_0)$ is even and $\alpha=1$   if $n(G_0)$ is odd, then  we  have 
\begin{eqnarray*}
\frac{|W_Q|}{c(Q-W_Q)} &  \le & 	\frac{f(Q_0)-1+2f_s+2n(Q_0)-b +\lfloor \frac{1}{2} n(Q_0) \rfloor }{e(Q_0)+3f_s+\lfloor \frac{1}{2} n(Q_0) \rfloor } \\
&=&   \frac{3n(Q_0)+0.5f_s+1-1.5 b-\alpha}{2n(Q_0)+1.5 f_s-0.5b-\alpha}. \\  
%	 &\le & \frac{3n(Q_0)+0.5f_s+1-0.75 b -0.25 n(Q_0)-\alpha}{2n(Q_0)+1.5 f_s-0.5b -\alpha} \quad (0.75b \ge 0.25n(Q_0))\\
%	 & \le& \frac{3n(Q_0)+1.5f_s-0.75 b -0.25 n(Q_0) -\alpha}{2n(Q_0)+1.5 f_s-\frac{b}{2}-\alpha} \\
%	 &=&  \frac{3n(Q_0)+2.25f_s-\frac{3}{4} b  -1.5\alpha-\frac{1}{4} n(Q_0) -0.75 f_s +0.5\alpha}{2n(Q_0)+1.5 f_s-\frac{b}{2}-\alpha}. 
\end{eqnarray*}
We claim that $b\ge 6$.  We already argued that $Q_0\cap C_0$ is connected. 
% If $Q_0\cap C_0$ has at least three components, then it is clear that $e_{G_0}(Q_0, G_0-V(Q_0)) \ge  e_{C_0}(Q_0\cap C_0, C_0-V(Q_0)) \ge 6$. 
%If $Q_0\cap C_0$ has  exactly two components, then we have $e_{C_0}(Q_0\cap C_0, C_0-V(Q_0)) \ge 4$ and $Q_0\cap C_0$
%has at least two vertices that are adjacent in $G_0$ to vertices from $\{u_1, \ldots, u_{34}\}$. Thus again we have $e_{G_0}(Q_0, G_0-V(Q_0)) \ge 6$.  
%Suppose now that $Q_0\cap C_0$ is connected. 
Then as any cycle of $G_0-w$ has length at least 8, we know that $Q_0$ contains at least 8 consecutive vertices of $C_0$. As $Q_0\cap C_0$ is connected and $Q_0$ is 2-connected, $Q_0$ contains a cycle that contains $Q_0\cap C_0$ as a subgraph.  As any cycle of $G_0-w$ that contains some 
consecutive vertices of $C_0$ has at least 6 vertices that are adjacent in $G_0$ to vertices from $\{u_1, \ldots, u_{34}\}$, we get 
$ e_{G_0}(Q_0, G_0-V(Q_0)) \ge 6$. 
 For any $u\in V(Q_0)$ such that $e_{G_0}(u, G_0-V(Q_0)) \ge 1$, we know that $R(u)$
intersects both $W$ and $Q$ by the same argument as in the first paragraph of this proof. Thus  $e_{G_0}(Q_0, G_0-V(Q_0)) \ge 6$ implies $b\ge 6$. 
Hence  $1.5(c(Q-W_Q)-1)-|W_Q| \ge 1.75f_s+0.75b-0.5\alpha-2.5 \ge 0$. 
%
%Since $n(Q_0) \ge 8$   by the fact   that  any 2-connected subgraph of $G_0-w$ contains at least 8 vertices, we have  $2n(Q_0)+1.5 f_s-\frac{b}{2}-\alpha \ge \frac{2}{3}(3n(Q_0)+2.25f_s-\frac{3}{4} b  -1.5\alpha-\frac{1}{4} n(Q_0) -0.75 f_s +0.5\alpha)+1$. 
Then   we have $|W\cup W_Q| >|W|$ and 
$h(W\cup W_Q)\ge h(W)$,   
giving a contradiction to the choice of $W$. 
\qed

By Claim~\ref{claim:one-2-connected-component}, $G-W$ has a  unique 2-connected component that contains vertices of $D$. We let $Q$ be  that  component.

\begin{CLA}\label{claim:maximum-number-of-components}
Let  $R$ be a subgraph of  $G-V(Q)$,  and  let $R_0$ be the component graph of $R$.  Then the vertices of $R$ can be separated in at most  $e(R_0)+ |V(R)\cap T_1|+\alpha'(R_0)$ components of $G-W$, where $\alpha'(R_0)$ is the size of a maximum matching in $R_0$.  
\end{CLA}

\pf  Every component of $R-W$  is either a single vertex or an edge  by Claims~\ref{claim:block-components} and~\ref{claim:one-2-connected-component}. 
Since each vertex from $S\cap V(R)$ has in $G$ (and also in $R$ by Claim~\ref{claim:outerface-component}) more than 
4 neighbors from $T\cup U$, Claims~\ref{claim:T-vertex-in-W} and~\ref{claim:C-vertex-in-W} imply that $V(R)\cap S \subseteq W$. 
 Hence  each component of 
$R-W$ is either a vertex from $T\cap V(R)$, an edge consisting  of one vertex 
from $T\cap V(R)$ and one vertex from $U_3\cap V(R)$,  or a vertex from $U_3\cap V(R)$. Note that  it is impossible to have a component of $R-W$ that is an edge consisting of two vertices  say $x$ and $y$, from $U_3$ by Claim~\ref{claim:T-vertices-adj-A},  as   if $xy$ were a component of $R-W$, then we must have $T(x), T(y) \in W$ where recall that $T(x)$ is the neighbor of $x$ from $T$ in $G$, but that gives a contradiction to Claim~\ref{claim:T-vertices-adj-A}. Also it is impossible to have a component in $R-W$ containing a vertex of $D$ as $D-W$ is contained in $Q$ by Claim~\ref{claim:W-no-intersect-B-C}. 

Each component of $R-W$ that contains a vertex from $T$ either corresponds to an edge of $R_0$ (if the vertex is from $T_2$) 
or is a vertex from $V(R)\cap T_1$. 
For each component of $R-W$ that is a single vertex $u$ from $U_3$, we have 
$T(u) \in W$.   Thus $T(u) \in T_2$ by Claim~\ref{claim:T-vertex-in-W}. Let $v$ be the other $U$-neighbor of $T(u)$.  Then both $u$ and $v$ are trivial components of $R-W$
by Claim~\ref{claim:T-vertices-adj-A}.  
As $uT(u)v$ corresponds to an edge of $R_0$,  the set of pairs of such components $u$ and $v$ of $R-W$ corresponds to a matching  of $R_0$.  Let $\alpha$ be the total number of pairs of such components $u$ and $v$.  
As the $\alpha$ vertices from $T$ that are adjacent  in $G$ to such pairs of vertices $u$ and $v$ are contained in $W$, 
 we know that $R-W$ has exactly $e(R_0)+|V(R)\cap T_1|-\alpha+2\alpha$ components.  The conclusion of 
 the claim now follows since $\alpha \le \alpha'(R_0)$. 
\qed

\begin{CLA}\label{claim:S-vertices}
Let $F$ be the boundary  of an $S$-triangle  $S^*$ with $S^*\subseteq W$. Then  vertices of $F-W$ are separated in at least 6 distinct components of $G-W$ and so $c(F-W) \ge 6$. 
\end{CLA}

\pf 
  We adopt the labeling of vertices from Figure~\ref{fig4} 
for convenient description.   Let $$F=w_1t_1x_{1,1}x_{1,2}t_2x_{2,1}x_{2,2}t_3x_{3,1} x_{3,2}t_4 x_{4,1} x_{4,2} t_5 x_{5,1} x_{5,2} t_6w_2w_1.$$
As $S^*\subseteq W$ and $Q$ is 2-connected, it follows 
that  none of  the three $T$-vertices respectively incident with $x_{1,3}, x_{3,3}, x_{5,3}$ 
is contained in $Q$. Furthermore, any of these $T$-vertices cannot be contained in  
$W$ by Claim~\ref{claim:T-vertices-adj-A}. Therefore, these three $T$-vertices are three components of $G-W$. As a consequence, we must have $x_{1,3}, x_{3,3}, x_{5,3} \in W$. 

Suppose otherwise that  vertices of $F-W$ are separated in at most 5 distinct components of $G-W$. 
Let  $|V(F)\cap W|  =j$  for some integer $j\ge 0$.  As $F$ is a cycle and so is 1-tough, we have $c(F-W)  \le j$, and so  vertices of $F$ are separated in at most $\max\{j,1\}$  distinct components of $G-W$ as well. 
Thus $j\le 5$ and 
vertices from  $V(F)\cap W\cup S^*\cup \{x_{1,3}, x_{3,3}, x_{5,3}\}$ are connected in $G$ to at most $\max\{1+3, j+3\}$  
components of $G-W$.   
We Claim that $W':=W\setminus (V(F)\cap W\cup S^*\cup \{x_{1,3}, x_{3,3}, x_{5,3}\})$ is still a cutset of $G$.  
It suffices to show that  $G-W$ has a component that is connected to vertices of $W'$ but is not connected to vertices of $V(F)\cap W\cup S^*\cup \{x_{1,3}, x_{3,3}, x_{5,3}\}$. 
Suppose instead that  vertices of $V(F)\cap W\cup S^*\cup \{x_{1,3}, x_{3,3}, x_{5,3}\}$ are connected to all the components of $G-W$.  If $j=0$, then we have $|W| \ge 6$ 
and $c(G-W) =4$ as vertices of $V(F)\cap W\cup S^*\cup \{x_{1,3}, x_{3,3}, x_{5,3}\}$ are connected to  exactly four  components of $G-W$. However this gives $h(W)=0$. 
Thus $j\ge 1$.  By Claims~\ref{claim:T-vertex-in-W} and~\ref{claim:C-vertex-in-W}, $W$ contains 
an $S$-vertex that is adjacent in $G$ to the vertices from $V(F)\cap W$ but is not contained in $S^*$. 
Thus $|W| \ge j+7$,   but we  get $h(W) \le \frac{3}{2} (j+3) -(j+7)  \le 0$ when $j\le 5$.  Therefore $G-W$ has a component that is connected to vertices of $W'$ but is not connected to vertices of $V(F)\cap W\cup S^*\cup \{x_{1,3}, x_{3,3}, x_{5,3}\}$, implying that 
$W'$ is a cutset of $G$. 
However,   as vertices from  $V(F)\cap W\cup S^*\cup \{x_{1,3}, x_{3,3}, x_{5,3}\}$ are connected in $G$ to at most $\max\{1+3, j+3\}$  
components of $G-W$  and $\max\{1+3, j+3\}<\frac{2}{3} (j+6)+1$ when $j\le 5$,  we get a contradiction to  Fact~\ref{fact:W-vertices-connect-to-components}.   
Thus vertices of $F-W$ are separated in at least 6 distinct components of $G-W$ and   so $c(F-W) \ge 6$.   
\qed

\begin{CLA}\label{claim:S-triangle-part-of-Q}
	Let $F$ be the boundary of an $S$-triangle  $S^*$ that is associated with $D$. Suppose that $S^*$  is contained in $W$,  then  $S^*_l, S^*_r, S^*_c\subseteq W$. 
	As a consequence, none of the spoke and the $C$-segment of $F$ is contained in $Q$. 
\end{CLA}

\pf We adopt the 
labels of vertices in Figure~\ref{fig4} for convenient description.   Thus  let  $$F=w_1t_1x_{1,1}x_{1,2}t_2x_{2,1}x_{2,2}t_3x_{3,1} x_{3,2}t_4 x_{4,1} x_{4,2} t_5 x_{5,1} x_{5,2} t_6w_2w_1.$$ 
By the same argument as in Claim~\ref{claim:S-vertices}, we have $x_{1,3}, x_{3,3}, x_{5,3} \in W$.  

If $S^*_c\not\subseteq W$, then we have $S_c^*\subseteq V(Q)$ by Claim~\ref{claim:outerface-component}. 
Then by the construction of of $G_0$ and Claim~\ref{claim:maximum-number-of-components}, we have $c(F-W) \le 5$, showing a 
contradiction to Claim~\ref{claim:S-vertices}. Thus $S^*_c\subseteq W$. 
We consider now,  by symmetry,  that $S^*_r\subseteq V(Q)$ and $S^*_l\subseteq W$.  
Suppose $w_2t_6x_{5,2}x_{5,1}t_5x_{4,2}$ is part of the boundary of $S^*_r$. By Claim~\ref{claim:S-vertices}, we need to have $c(F-W) \ge 6$. 
This in particular, implies that $x_{4,1} \in W$. 
 Let $W_Q=\{w_2, x_{5,2}, t_5, x_{4,3}, x\} \cup S^*_r$, where  $x$ is the other $U$-neighbor of $T(x_{4,3})$. 
 Then $Q-W_Q$ has at least 5 components: the vertex $t_6$, the vertex $x_{5,1}$, the vertex $x_{4,2}$, the vertex $T(x_{4,3})$, 
 and the component containing a vertex of $D$. 
Since  $|W_Q|=6$ and $c(G-(W\cup W_Q))  \ge c(G-W)+4$,  and so $h(W\cup W_Q) \ge h(W)$,  we get a contradiction to the choice of $W$. 
Therefore $S^*_l, S^*_r, S^*_c\subseteq W$.  

We argue that  it is impossible to have any  of the two spokes contained in $F$ to be contained in $Q$. Suppose, by symmetry, that $w_2t_6x_{5,2}x_{5,1}t_5x_{4,2}$  is contained in $Q$. 
Then as  $c(F-W) \ge 6$, we must have $x_{4,1} \in W$.  Since $Q$ is 2-connected and $S^*_l, S^*_r, S^*_c\subseteq W$, it follows that 
the $U$-neighbor, say $x$ of $T(x_{4,3})$ is contained in $Q$. Let $W_Q=\{w_2, x_{5,2}, t_5, x_{4,3}, x\}$. 
Then again we get  $h(W\cup W_Q) \ge h(W)$, a contradiction to the choice of $W$.  
Lastly, the $C$-segment $x_{2,1} x_{2,2}t_3x_{3,1}x_{3,2} t_4x_{4,1}x_{4,2}$  of $F$ is not contained in $Q$ 
since otherwise $c(F-W)  \le 5$ by Claim~\ref{claim:maximum-number-of-components}. 
\qed

\begin{CLA}\label{claim:S-triangle-part-of-Q2}
	Let $F$ be the boundary of an $S$-triangle   $S^*$  that is associated with $D$. Suppose $S^*\subseteq V(Q)$, $S^*$ is internal,  and one of $S_l^*$  or $S_r^*$  is contained in $W$. 
	 Then $S_c^*\subseteq V(Q)$. 
\end{CLA}

\pf  We adopt the 
labels of vertices in Figure~\ref{fig4} for convenient description.   Thus  let  $$F=w_1t_1x_{1,1}x_{1,2}t_2x_{2,1}x_{2,2}t_3x_{3,1} x_{3,2}t_4 x_{4,1} x_{4,2} t_5 x_{5,1} x_{5,2} t_6w_2w_1.$$ 
Since $S^*\subseteq V(Q)$, we have $F\subseteq Q$ and $x_{1,3}, x_{3,3}, x_{5,3}, T(x_{1,3}), T(x_{3,3}), T(x_{5,3}) \in V(Q)$   by Claims~\ref{claim:T-vertex-in-W} and~\ref{claim:C-vertex-in-W}.

Suppose, by symmetry,  that  $S^*_l\subseteq W$ and  $w_2t_6x_{5,2}x_{5,1}t_5x_{4,2}$ is part of the boundary of $S^*_r$.  Suppose to the contrary that $S_c^*\subseteq W$. 
Consider first that  $x_{2,3} \in W$.  
Let $W_Q=\{x_{1,3}, x_{3,3}, x_{5,3},  w_1, x_{1,1}, t_2, x_{2,2}, x_{3,1}, x_{4,1}\} \cup S^*$. 
Then we know that $|W_Q|=12$ and $c(G-(W\cup W_Q))  \ge c(G-W)+8$ and so $h(W\cup W_Q) \ge h(W)$, a contradiction to the choice of $W$. 

Thus we suppose $x_{2,3} \not\in  W$. 
Let $F_l$ be the boundary of $S_l^*$. 
 Then $x_{2,3} \in V(Q)$. As a consequence, we have $t_7 \in V(Q)$ by Claim~\ref{claim:T-vertices-adj-A} and so $z_{1,2} \in V(Q)$ by Claim~\ref{claim:C-vertex-in-W}. 
We claim that $z_{1,1} \in W$. For otherwise, by the same logic as in the line above, we have $t_8, y_{4,2} \in V(Q)$. 
Let $S'$ be the $S$-triangle with boundary $F'=w_3t_9y_{1,1}y_{1,2}\ldots y_{5,1}y_{5,2} t_{14}w_4w_3$.  
Then we must have $S'\in W$: otherwise $F' \subseteq Q$, and so vertices of $S_l^*$ are only connected in $G$ to the component $Q$ of $G-W$, contradicting Fact~\ref{fact:W-vertices-connect-to-components}. 
Now as   $S'$ is in $W$, we have $S_l', S_r', S_c'\in W$ by Claim~\ref{claim:S-triangle-part-of-Q}, and none of the spokes of $F'$ is contained in $Q$ and the $C$-segment of $F'$ is not  contained in $Q$ as well. 
Since none of the spokes of $F'$ is contained in $Q$, we have $w_3, w_4\in W$ by Claim~\ref{claim:spoke}.  Since the $C$-segment of $F'$ is not contained  in $Q$, some vertex of the $C$-segment of $F'$ is contained in $W$. 
Now as  $S_l^* S_l', S_r', S_c'\in W$, it follows that the vertex $z_{1,1}$ is a cutvertex of $Q$ (deleting $z_{1,1}$ separates $t_8$ with $z_{1,2}$), a contradiction to $Q$ being 2-connected.  

Thus we have $z_{1,1} \in W$. By the same argument as above, we must have $S'\in W$: otherwise $F' \subseteq Q$, and so the maximum number of components from 
the boundary of $S_l^*$ we can  get is by deleting $y_{4,3}$ and $z_{1,1}$ to get $t_8$ and the rest as two components.   As $G$ is  3-connected and $h(W)>0$, we know that $c(G-W)\ge 3$. 
As  we deleted the vertex in $S_l^*$ and  at most two  vertices $y_{4,3}$ and $z_{1,1}$ from $F_l$ and these three vertices 
are connected in $G$ to at most two distinct components of $G-W$, we know that $W\setminus (S_l^* \cup \{y_{4,3}, z_{1,1}\})$ is a cutset of $G$. However, 
 $2<\frac{2}{3}\times 3+1$ contradicts Fact~\ref{fact:W-vertices-connect-to-components}.  Thus    $S'$ is in $W$, and so  $S_l', S_r', S_c'\in W$ by Claim~\ref{claim:S-triangle-part-of-Q}. 

Let   $u$ be the  other $U$-neighbor of $T(z_{1,3})$, and $$W_Q=\{z_{1,2}, u, x_{2,2}, x_{2,3}, x_{3,1}, x_{3,3}, x_{4,1}, t_2, x_{1,1}, x_{1,3}, x_{5,3}, w_1\} \cup S^*.$$
Then we created 10 extra components after deleting all the 10 vertices of $W_Q$. Those  components include single vertex components and components consisting of 
an  edge:  $$t_1, x_{1,2}, x_{2,1}, t_7, T(z_{1,3})z_{1,3}, t_3, x_{3,2}t_4$$ and the three vertices from $T_1$ that are  associated with $S^*$. 
However we get  $h(W\cup W_Q) \ge h(W)$, a contradiction to the choice of $W$. 
\qed

Let $S^*$ be an $S$-triangle associated with $D$. If $S^* \subseteq V(Q)$, then the two spokes  of $G$ that are contained in the boundary  of $S^*$ are both contained in $Q$ by 
Claims~\ref{claim:T-vertex-in-W} and~\ref{claim:C-vertex-in-W}. If $S^* \not\subseteq V(Q)$ then we have $S^*\subseteq W$ by Claim~\ref{claim:outerface-component}. 
By Claim~\ref{claim:S-triangle-part-of-Q},   none of the two spokes contained in the boundary of  $S^*$  is contained in $Q$, and the $C$-segment of the boundary of $S^*$
is also not contained in $Q$. 

We say that two  $S$-triangles $S^*$ and $S'$  associated with $D$ are in the same \emph{patch} if $S_c^*=S_c'$. A \emph{patch of $S$-triangles} is the union of all the $S$-triangles that are from the same patch. 
Similarly, a \emph{patch of spokes} is the union of the spokes contained in the boundaries of a patch of $S$-triangles. 
By Claim~\ref{claim:S-triangle-part-of-Q2} and the argument immediately above, a patch of $S$-triangles are either all contained in $Q$ or all contained in $W$.  In particular, 
when a patch of $S$-triangles are all contained in $Q$, then by Claims~\ref{claim:T-vertex-in-W} and~\ref{claim:C-vertex-in-W}, all the spokes that are contained in the boundary of this  patch of $S$-triangles are contained in $Q$ as well.  When a patch of $S$-triangles are all contained in $W$, none of the spokes from the  corresponding patch is 
contained in $Q$ by Claim~\ref{claim:S-triangle-part-of-Q}. 

Let  $Q_0$ be the component graph of $Q$, $a=|V(Q)\cap V(A)|$, $b$ be the number of $\CC_3$-triangles that intersects both $W$ and $Q$, and $c$ be the number of components of $Q\cap C$.  
As $Q$ is 2-connected, Claim~\ref{claim:outerface-component} implies that $Q_0$ is 2-connected. 
We first claim that $b\ge 2c$.

\begin{CLA}\label{b>=2}
It holds that $b\ge 2c\ge 2$. 
\end{CLA}

\pf   The conclusion is obvious if $Q_0\cap C_0  \ne C_0$, as  in this case we have $b=e_{G_0}(Q_0-w, W\cap V(G_0)) \ge e_{C_0}(Q_0\cap C_0, W\cap V(G_0))  \ge 2c$. 
Thus we may assume that $Q_0\cap C_0  = C_0$ and $b\le 1$.   This implies that at most one spoke of $G$ is not contained in $Q$. 
By Claim~\ref{claim:S-triangle-part-of-Q}, every $S$-triangle associated with $D$ is contained in $Q$, and the same holds for all the $S$-triangle not associated with $D$ as a patch of $S$-triangles are either all contained in $Q$ or all contained in $W$.  Thus 
 $T_1\subseteq V(Q)$ by Claim~\ref{claim:T-vertex-in-W}.  
 %Let $W_0$  be the vertex set of the component graph of $G[W]$. 
Then as $Q_0\cap C_0  = C_0$,  it follows that $G_0-V(Q_0) -w$  consists of isolated vertices. This together with the fact that every component of $G-W-V(Q)$ 
that is a single vertex from $U_3$ corresponds to an edge of $G_0-V(Q_0) -w$,  it implies  that  no component of $G-W-V(Q)$ is a single vertex from $U_3$. 
Thus every  component  of $G-W-V(Q)$
is either a vertex from  $T_2$ or an edge consisting of a vertex from $T_2$ and a vertex from $U_3$. 
Since $b\le 1$ and  $Q_0\cap C_0  = C_0$, it follows that $Q_0$ contains all the edges that incident in $G_0$ with a vertex of $C_0$. 
Then as  at most one spoke of $G$ is not contained in $Q$, it follows that $c(G-W) \le 3$. 
As $G$ is 3-connected and $h(W)>0$, it follows that $c(G-W)=3$. 
Thus  exactly one  long spoke, say $P$, of $G$ is not contained in $Q$.   Since  we need to delete the two endvertices of  $P$  for it to be not 
contained in $Q$ by Claim~\ref{claim:spoke},  at least three vertices  of $P$  need  be deleted  to  disconnect  $P$ into two components. 
We also need to delete at least two $S$-vertices that have $P$ has part of their boundaries. Thus $|W| \ge 5$, showing a contradiction to $h(W)>0$. 
\qed

In the next, we will show that we can construct a cutset $W_Q$ of $Q$ that contains $V(A)\cap V(Q)$ and that $1.5(c(Q-W_Q)-1) -|W_Q| \ge -0.75$. As $h(W)>0$ implies $h(W) \ge 0.5$, 
it will  then follow  that $h(W\cup W_Q) \ge -0.25$.  This will lead to a contradiction as $h(W\cup W_Q)$ should be at most $-0.5$.  
If $a=0$, then we let $W_Q=\emptyset$, which certainly satisfying $1.5c(Q-W_Q) -|W_Q| \ge 0.75$. Thus we assume $a\ge 1$
and construct the set $W_Q$. Let $f_s$ be the number of $S$-triangles contained in $Q$. 

By the same argument as in the proof of Claim~\ref{claim:one-2-connected-component} and Euler's formula, we have 
\begin{eqnarray*}
e(Q_0)  & =&  0.5(3(n(Q_0)-1)+a-3f_s-b)=1.5 n(Q_0)+0.5 a-1.5f_s-0.5b-1.5, \\
f(Q_0)&=&e(Q_0)-n(Q_0)+2 = 0.5n(Q_0)+0.5 a-1.5f_s-0.5b+0.5. 
\end{eqnarray*}

Let $M_0$ be a maximum matching in  $Q_0$, and let $a'$ be the number of short spokes that are contained in $Q$.   
As $c(Q_0\cap C_0)=c(Q\cap C)$ and each component of $Q_0\cap C_0$  is a path, we know that  $|M_0| \ge 0.5(n(Q_0)-(a-a')-c)$. 
For each $uv\in M_0$,  there exists $x_u\in R(u)$ and $x_v\in R(v)$
such that $x_u$ and $x_v$ are  both adjacent in $G$ 
to a vertex $y\in T$.  We call $x_u$  a \emph{representative vertex}
of $R(u)$.  As $M_0$ is a matching, each  graph  from 
$\CC_3$   either has no representative vertex or has 
a unique representative vertex. 

Let $R\in \CC_{3}$ such that $V(R)\cap V(Q) \ne \emptyset$.  By the argument as in the first paragraph in the proof of Claim~\ref{claim:one-2-connected-component}, we have that either $R\subseteq Q$ or 
$|V(R) \cap V(Q)| =2$. 
If $R$ has a representative vertex, say  $x$,  
let $W_R \subseteq (V(R) \cap V(Q))\setminus \{x\}$  be the  set of 
two   vertices (if  $R\subseteq Q$) or one  vertex (if $|V(R) \cap V(Q)|=2$) such that $e_G(R-W_R, T)=e_G(x,T)$. 
Otherwise,  let $W_R \subseteq V(R)$  be a  set of 
two   vertices (if  $R\subseteq Q$) or one  vertex (if $|V(R) \cap V(Q)|=2$) such that  $e_G(R-W_R, T)=1$. 
Note that for any two representative vertices $x_1$ and $x_2$,  $x_1$
and $x_2$ are adjacent in $G$ to the same vertex 
from $T$. We let $T^*$ be the set of all 
these vertices from $T$ that are adjacent in $G$
to a representative vertex of  components from $\CC_3$ that intersects $Q$.

%, and so $f(Q_0)=e(Q_0)-n(Q_0)+2 =0.5n(Q_0)-1.5f_0-\frac{1}{2}b+2$.  

Let $S^*$ be the set of  $S$-vertices that are embedded inside a face with boundary, say $F$.  By Claim~\ref{claim:outerface-component}, we have either $S^*\subseteq V(Q)$ or $S^*\cap V(Q) =\emptyset$. If $S^*\subseteq V(Q)$, then we have $F\subseteq Q$ by Claims~\ref{claim:T-vertex-in-W} and~\ref{claim:C-vertex-in-W}.  Thus $Q_0$ has a face 
whose boundary is the component graph of $F$.  Therefore we have  $|V(Q)\cap S| \le f(Q_0)-1+2f_s$,  as at least one face of $Q$ whose boundary has vertices adjacent in $G$ to vertices from $W$ and so 
there is no   $S$-vertex  embedded in $Q$ inside that face.

%For any $D\in \CC_{2k+1}^1$ for some $k\ge 2$, by Claim~\ref{claim:large-odd-comp},  we let $S_D \subseteq V(D)$  be a  set of  $\lfloor \frac{4k+2}{3} \rfloor$ vertices
% such that $D-S_D$ is set of isolated vertices. 

%For any $D\in \CC_{2k+1}^2$ for some $k\ge 2$, we let $S_D \subseteq V(D)$  be the  set of 
%$2k+1$ vertices such that $e_G(D-S_D, T)=0$. 
%Note that $D-S_D$ is a graph with at least one vertex. 

%As there are $2m$ representative vertices, we have $|T^*| = m$. 
%Then we have $|M_0| \ge \frac{1}{2}n(Q_0)$. 
Let $W_Q$ be the set that consists of all vertices in  $V(A)\cap Q$, $S\cap V(Q)$,  $T^*$, and $W_R$ for all $R\in \CC_3$ such that $|V(R)\cap V(Q)| \ge 2$.  
Then 
\begin{eqnarray*}
|W_Q|  &\le& (2(n(Q_0)-1)-b)+a+f(Q_0)-1+2f_s+|T^*|\\
 &=& 2n(Q_0)+a+f(Q_0)+2f_s+|T^*|-b-3 \\
 &=& 2.5n(Q_0)+1.5a+0.5f_s+|T^*|-1.5b-2.5. 
\end{eqnarray*}

In $Q-W_Q$, 
there are exactly $e(Q_0)-|T^*|$ components that contains a vertex of $T_2$, $3f_s$ components that contains a vertex of $T_1$, and  $2|T^*|$ components that each is a single vertex from $U_3$, 
 and $D-W-W_Q$ is also a component. Thus we have 
 \begin{eqnarray*}
 	c(Q-W_Q) &=&(e(Q_0)-|T^*|)+ 3f_s +2|T^*|+1\\
 	&=& 1.5n(Q_0)+0.5a+1.5f_s  +|T^*|-0.5b-0.5. 
 \end{eqnarray*}
%since  the three vertices from $T$ that are adjacent in $G$ to vertices from every $S$-triangle of $Q$ are contained in $Q$ by Claim~\ref{claim:T-vertex-in-W}, and $D-W-W_Q$ is also a component.  
% As $n(Q_0) \ge 3f_s+b$, it then follows that $2.5n(Q_0)+1.5a+0.5f_s -1.5 b -2.5 \ge 1.5n(Q_0)+0.5a+1.5f_s-0.5b  -0.5$. 
 Suppose that $Q$ contains precisely $h$ patches of $S$-triangles, where $h\in [0,4]$. 
Then we have $a=2(f_s-h)+a'$, as $Q$ contains $h$ $S$-triangles that are not associated with $D$. 
Then  as $|T^*| \ge \frac{1}{2}(n(Q_0)+a'-a-c)$,  we get 
\begin{eqnarray*}
	&&1.5c(Q-W_Q)-|W_Q| \\
%	&=& 1.5\left(e(Q_0)+3f_s+|T^*|+1\right) -\left(2n(Q_0)+a+f(Q_0)+2f_s+|T^*|-b-3\right) \\
	&=& 1.5(1.5n(Q_0)+0.5a+1.5f_s  +|T^*|-0.5b-0.5)\\
	&&-(2.5n(Q_0)+1.5a+0.5f_s+|T^*|-1.5b-2.5) \\
	&=& -0.25n(Q_0)-0.75a+1.75f_s+0.5|T^*|+0.75b+1.75 \\
	&\ge & -0.25n(Q_0)+1.75f_s-0.75a+0.25(n(Q_0)+a'-a-c)+0.75b+1.75 \\
  &=& 1.75f_s+0.25a'+0.75b+1.75-a-0.25c \\ 
  &\ge &-0.25f_s-0.75a'+2h +0.625b+1.75   \quad   \text{($a=2(f_s-h)+a'$ and $b\ge 2c$)}.  
\end{eqnarray*}

%Since $Q\ne Q_0$ and $Q_0$ is 2-connected, it follows that $b\ge 2$. 
Note that $f_s\le 21$, $b\ge 2c \ge 2$, and $a'\le 5$. 
When $h=4$,  we have 
\begin{eqnarray*}
-0.25f_s-0.75a'+2h +0.625b+1.75 &\ge& -5.25-3.75+8+1.25+1.75 \ge 2. 
\end{eqnarray*} 
Thus $1.5(c(Q-W_Q)-1)-|W_Q| >0$ and so $h(W\cup W_Q)>h(W)$, a contradiction to the choice of $W$. 
Thus we have $h\le 3$. 
 
 Suppose $h=3$. Then we have 
 $f_s\le 17$ as there are  at least four $S$-triangles in a patch.  Thus 
 \begin{eqnarray*}
 	-0.25f_s-0.75a'+2h +0.625b+1.75 \ge -4.25-3.75+6+1.25+1.75 \ge 1. 
 \end{eqnarray*}
%Thus $-0.25f_s-0.75a'+2h +0.625b+1.75 \ge -4.5-3.75+6+1.25+1.75 \ge 0.75$. 
Thus $1.5(c(Q-W_Q)-1)-|W_Q|  \ge -0.75$. 

% Since $h(W) \ge 0.5$, it follows that $h(W\cup W_Q) \ge -0.25$. 
%This is a contradiction since we in fact  have $h(W\cup W_Q)  \le -0.5$.  

Next we suppose $h=2$.   
Then we have $f_s\le 13$. If $a'\le 4$, then  we get $1.5(c(Q-W_Q)-1)-|W_Q|  \ge -0.75$. 
Thus we suppose $a'=5$.  Since two patches of $S$-triangles are not contained in $Q$ but $Q$ contains all of the 5 short spokes, 
it follows that $c\ge 2$. Thus $b\ge 2c\ge 4$. Hence 
\begin{eqnarray*}
-0.25f_s-0.75a'+2h +0.625b+1.75 &\ge& -3.25-3.75+4+2.5+1.75 =1.25. 
\end{eqnarray*}
Thus again we get $1.5(c(Q-W_Q)-1)-|W_Q|  \ge -0.25 \ge -0.75$.

Then  we suppose $h=1$.  We claim that $f_s\le 5+c$.     A largest component of $Q_0\cap C_0$ can contain the $C$-segments of  boundaries of at most six  $S$-triangles (the component that 
contain the $C$-segment of the boundary of an $S$-triangle from the unique patch that is contained in $Q$ ). Every other component of $Q_0\cap C_0$ does not contain any $C$-segment of the boundary of any  $S$-triangle from a patch 
and so it 
 can contain 
 the $C$-segment  of   the boundary of at most one  $S$-triangle.  Hence we have  $f_s\le 5+c$. 
 Since $c\le 0.5b$,  when $a'\le 3$, we get 
 \begin{eqnarray*}
 -0.25f_s-0.75a'+2h +0.625b+1.75  & \ge&  -1.25-0.75a'+0.5b+3.75  \ge 1.25. 
 \end{eqnarray*} 
Thus $a'\ge 4$. This implies that $b\ge 2c\ge 4$. Then again, $-1.25-0.75a'+0.5b+3.75 \ge 0.75$.

Lastly, suppose $h=0$. Then we have $a\le 2f_s+c+1$ and $f_s\le c+1$  by the construction of $G_0$. 
Then  if $a'\ge 1$,  as $b\ge 2c$, we get 
\begin{eqnarray*}
1.5c(Q-W_Q)-|W_Q| &\ge & 1.75f_s+0.25a'+0.75b+1.75-a-0.25c \\
 &\ge & -0.25f_s-1.25c+0.25a'+0.75b+0.75 \\
 &\ge & -1.5c+0.25a'+0.75b+0.5  \ge 0.25a'+0.5 \ge 0.75. 
\end{eqnarray*}

Thus we suppose $a'=0$.  As a consequence, we have $a=2f_s$. Then 
\begin{eqnarray*}
	1.5c(Q-W_Q)-|W_Q| &\ge & 1.75f_s+0.25a'+0.75b+1.75-a-0.25c \\
	&\ge & -0.25f_s+0.25a'+0.75b+1.75-0.25c \\
	&\ge & -0.25f_s+0.625b+1.75 \ge 1.75, 
\end{eqnarray*}
as $f_s\le 5$. Thus $1.5(c(Q-W_Q)-1)-|W_Q| >0$ and so $h(W\cup W_Q)>h(W)$, a contradiction to the choice of $W$.

By the arguments above, we can always find $W_Q$ with $1.5(c(Q-W_Q)-1) -|W_Q|  \ge -0.75$. 
Let $W^*=W\cup W_Q$.  Then as $h(W)>0$ implies $h(W) \ge 0.5$ by the definition of the function $h$, it follows that $h(W^*) \ge -0.25$. 
Now we achieve a contradiction   by showing  that $h(W^*) \le -0.5$. 

Note that by  the assumption that $Q$ is the only 2-connected component of $G-W$, we know that 
$D-W^*$ is the only 2-connected component of $G-W^*$ and $ V(A) \subseteq W^*$. 
Let $T^*=T\cap W$. Then by the same argument as above for calculating $|W_Q|$ and $c(Q-W_Q)$, and noting that $S\subseteqq W$, we have 
\begin{eqnarray*}
	|W^*|  &=& 2(n(G_0)-1)+ |V(A)| +f(G_0)+2f_s(G_0)+|T^*|\\
	&=& 240+39+49+42+|T^*|
\end{eqnarray*}
and 
$$c(G-W^*) =(e(G_0)-|T^*|)+3f_s(G_0)+2|T^*|+1 =168+63+|T^*|+1.$$
As $|T^*| \le 43$ by Fact~\ref{fact:matching-size-in-G0}, 
we get 
\begin{eqnarray*}
h(W^*)=1.5c(G-W^*)-|W^*|  & =& 1.5(232+|T^*|)-(370+|T^*|) \\
 &= & -22+0.5|T^*|  \le -22+21.5  \le -0.5, 
\end{eqnarray*}
giving a contradiction. 
The proof of Theorem~\ref{main} is now completed. 
\qed 

\section*{Acknowledgments}
The author is very grateful to the two anonymous referees  for their careful reading and valuable
comments.  The author 
also thanks Drs. Mark Ellingham, Dong Ye and Xiaoya Zha  for helpful discussions 
in the early stage of this project.   The author was supported by National Science Foundation grant 
DMS-2345869.

%\bibliographystyle{plain}
%\bibliography{SSL-BIB}

\end{document}